\documentclass[aos,MSNbibl,seceqn,dvips]{arximspdf}
\usepackage{graphicx}

\doi{10.1214/12-AOS1029} %kopijuoti is PTS
\volume{40}
\issue{4}
\pubyear{2012}
\firstpage{2069}
\lastpage{2101}

\makeatletter
\newcommand{\rrVert}{\Vert}
\newcommand{\llVert}{\Vert}

\newtheorem{theorem}{Theorem}[section]
\newtheorem{lemma}{Lemma}[section]
\newtheorem{corollary}{Corollary}[section]
\newtheorem{proposition}{Proposition}[section]

\newproclaim{example}{Example}[section]

\newproclaim{remark}[theorem]{Remark}

\newcommand{\given}{|}
\newcommand{\ra}{\rightarrow}
\newcommand{\rn}{\sqrt{n}}
\newcommand{\E}{\mathrm{E}}
\newcommand{\al}{\alpha}
\newcommand{\be}{\beta}
\renewcommand{\b}{\beta}
\newcommand{\ga}{\gamma}
\newcommand{\la}{\lambda}
\renewcommand{\l}{\lambda}
\newcommand{\epl}{\underline{\veps}_n}
\newcommand{\h}{\eta}
\renewcommand{\k}{\kappa}
\newcommand{\s}{\sigma}
\newcommand{\te}{\theta}
\newcommand{\q}{\theta}
\newcommand{\m}{\mu}
\newcommand{\veps}{\varepsilon}
\newcommand{\pli}{+\infty}
\newcommand{\A}{\cA}
\newcommand{\cA}{{\mathcal{A}}}
\newcommand{\cS}{{\mathcal{S}}}

\newcommand{\RR}{\mathbb{R}}
\newcommand{\NN}{\mathbb{N}}

\makeatother

\begin{document}
\begin{frontmatter}

\title{Needles and Straw in a Haystack:
Posterior concentration for possibly sparse sequences\thanksref{T1}}
\runtitle{Sparsity and Bayes posterior measure}

\thankstext{T1}{Supported in part by a Postdoctoral fellowship
from the VU University Amsterdam and
ANR Grant ``Banhdits'' ANR-2010-BLAN-0113-03.}

\begin{aug}
\author[A]{\fnms{Isma\"el} \snm{Castillo}\corref{}\ead[label=e1]{ismael.castillo@upmc.fr}}
\and
\author[B]{\fnms{Aad} \snm{van der Vaart}\ead[label=e2]{aad@cs.vu.nl}}
\runauthor{I. Castillo and A. W. van der Vaart}
\affiliation{Universit\'es Paris VI \& VII and VU University}
\address[A]{CNRS \& Laboratoire de Probabilit\'es\\
\quad et Mod\`eles Al\'eatoires\\
Universit\'es Paris VI \& VII\\
175, rue du Chevaleret\\
Paris, 75013\\
France\\
\printead{e1}}
\address[B]{Department of Mathematics\\
Faculty of Sciences\\
VU University\\
Amsterdam\\
The Netherlands\\
\printead{e2}} %adresu isvedimo komanda gale!
\end{aug}

% HISTORY:
\received{\smonth{4} \syear{2011}}
\revised{\smonth{11} \syear{2011}}

% ABSTRACT
%
\begin{abstract}
We consider full Bayesian inference in the multivariate normal mean
model in the situation that the mean vector is sparse. The prior
distribution on the vector of means is constructed hierarchically by
first choosing a collection of nonzero means and next a prior on the
nonzero values. We consider the posterior distribution in the
frequentist set-up that the observations are generated according to a
fixed mean vector, and are interested in the posterior distribution of
the number of nonzero components and the contraction of the posterior
distribution to the true mean vector. We find various combinations of
priors on the number of nonzero coefficients and on these coefficients
that give desirable performance. We also find priors that give
suboptimal convergence, for instance, Gaussian priors on the nonzero
coefficients. We illustrate the results by simulations.
\end{abstract}

% KEYWORDS
%
\begin{keyword}[class=AMS]
\kwd{62G05}
\kwd{62G20}
\end{keyword}
\begin{keyword}
\kwd{Bayesian estimators}
\kwd{sparsity}
\kwd{Gaussian sequence model}
\kwd{mixture priors}
\kwd{asymptotics}
\kwd{contraction}
\end{keyword}

\end{frontmatter}

%s1 #&#
\section{Introduction}
Suppose that we observe a vector $X=(X_1,\ldots,X_n)$
in $\RR^n$ such that
%
%e1.1 #&#
%
\begin{equation}
\label{EqModel} X_i = \theta_i +
\veps_i,\qquad i=1,\ldots,n,
\end{equation}
for independent standard normal random variables $\veps_i$ and
an unknown vector of means $\theta=(\te_1,\ldots,\te_n)$. We are
interested in Bayesian inference on $\theta$,
in the situation that this vector is possibly \textit{sparse}.

Non-Bayesian approaches to this problem have recently been considered
by many authors. Golubev~\cite{Golubev02} obtained results for model
selection methods and threshold estimators for the mean-squared
risk. Birg\'e and Massart~\cite{BirgeMas01} treated the model within
their general context of model selection by penalized least
squares. Abramovich et al.~\cite{fdr06} studied the performance of
the False Discovery Rate method. The\vadjust{\goodbreak} earlier work by Donoho and
Johnstone~\cite{DJ94} can be viewed as studying the problem within an
$\ell_r$ context. Many authors
(see, e.g.,~\cite{BRT09,zou,zhanghuang} and
references cited
there) have investigated the connection to the LASSO or similar methods.

Methods with a Bayesian connection were studied by George and Foster~\cite{GeorgeFoster}, Zhang~\cite{zhang}, Johnstone and Silverman
\cite{JohnSilv04,JohnSilvWave05}, Abramovich, Grinshtein and Pensky
\cite{abramovich2} and Jiang and Zhang~\cite{jiangzhang}. George and Foster
\cite{GeorgeFoster} and Johnstone and Silverman~\cite{JohnSilv04}
considered an empirical
Bayes method, consisting of modeling the parameters $\q_1,\ldots,\q_n$
a priori as independently drawn from a mixture of a Dirac
measure at $0$ and a continuous distribution, determining an
appropriate mixing weight by the method of (restricted) marginal
maximum likelihood and finally employing the posterior median or
mean. The second paper~\cite{abramovich2} motivated penalties, applied
in a
penalized minimum contrast scheme, by prior distributions on the
parameters, and derived estimators for the number of nonzero $\q_i$
and the $\q_i$, itself. The first is a posterior mode, but the
estimator for $\q$, called ``Bayesian testimation,'' does not seem
itself Bayesian. (In fact, the Gaussian prior for the nonzero
parameters in~\cite{abramovich2} will be seen to perform suboptimally
in our
fully Bayesian set-up.) Zhang~\cite{zhang} and
Jiang and Zhang~\cite{jiangzhang} obtain sharp results
on (nonparametric) empirical Bayes estimators.

Other related papers include \cite
{yuanlin,caijinlow,candestao,HuangMaZhang,jiangzhang,browngreenshtein}.
%,~\cite{abramovich2}.

A penalized minimum contrast estimator can often be viewed as the mode
of the posterior distribution, and it is helpful to interpret
penalties accordingly. However, the Bayesian approach yields a full
posterior distribution, which is a random probability distribution on
the parameter space. It has both a location and a spread, and can be
marginalized to give posterior distributions for any functions of the
parameter vector of interest. It is this object that we study in this
paper. Such full Bayesian inference was recently considered by Scott
and Berger~\cite{ScottBerger10}, who discussed various aspects not
covered in the present paper, but no concentration results. One
example of our results is that the beta-binomial priors in
\cite{ScottBerger10}, combined with moderately to heavy tailed priors
on the nonzero means, yield optimal recovery.

\textit{Sparsity} can be defined in various ways. Perhaps the most natural
definition is the class of \textit{nearly black} vectors, defined as
\[
\ell_0[p_n] = \bigl\{ \theta\in\RR^n\dvtx
\#(1\le i\le n\dvtx\te_i \neq0 ) \leq p_n \bigr\}.
\]
Here $p_n$ is a given number, which in theoretical investigations
is typically assumed to be $o(n)$, as $n\ra\infty$. Sparsity may also
mean that many means are small, but possibly not exactly zero.
Definitions that make this precise use
\textit{strong} or \textit{weak $\ell_s$-balls}, typically for $s\in(0,2)$.
These are defined as, with
$\te_{[1]}\geq\te_{[2]}\geq\cdots\geq\te_{[n]}$ the nonincreasing
permutation of the coordinates of $\theta=(\te_1,\ldots,\te_n)$,
\begin{eqnarray*}
\ell_s[p_n] & = & \Biggl\{ \theta\in\RR^n\dvtx
\frac1n\sum_{i=1}^{n} |
\theta_i|^s \leq\biggl(\frac{p_n}n
\biggr)^s \Biggr\},
\\
m_s[p_n] & = & \biggl\{ \theta\in\RR^n\dvtx
\frac{1}n\max_{1\le i\le n} i|\te_{[i]}|^s \leq
\biggl(\frac{p_n}n \biggr)^s \biggr\}.
\end{eqnarray*}
Because the nonzero coefficients in $\ell_0[p_n]$ are not
quantitatively restricted, there is no inclusion relationship
between this space and the weak and strong balls,
although results for the latter can be obtained
by projecting them into $\ell_0[p_n]$.
On the other hand,
the inclusion $\ell_s[p_n] \subset m_s[p_n]$ holds for any $s>0$.

The extent of the sparsity, measured by the constant $p_n$, is assumed
unknown. Our Bayesian approach starts by putting a prior $\pi_n$ on
this number,
a given probability measure on the set $\{0,1,2,\ldots, n\}$.
Next we complete this to a prior on the set of
all possible sequences $\theta=(\te_1,\ldots,\te_n)$ in $\RR^n$, by
given a draw $p$ from~$\pi_n$, choosing a random subset $S\subset
\{1,\ldots,n\}$ of cardinality $p$, and choosing the
corresponding coordinates $(\te_i\dvtx i \in S)$ from a density $g_S$ on
$\RR^S$
and setting the remaining coordinates $(\te_i\dvtx i\in S^c)$ equal to zero.
Given this prior, Bayes's rule yields the posterior distribution of
$\theta$, as usual. We investigate the properties of this posterior
distribution, in its dependence on the priors on the dimension and
on the nonzero coefficients, in the non-Bayesian set-up
where $X$ follows~(\ref{EqModel}) with $\q$ equal to a fixed,
``true'' parameter $\q_0$.

If the true parameter vector $\q_0$ belongs to $\ell_0[p_n]$, then it
is desirable that the posterior distribution concentrates most of its
mass on nearly black vectors. One main result of the paper is that
this is the case provided the prior probabilities $\pi_n\{p\}$
decrease exponentially fast with the dimension $p$.

The quality of the reconstruction of the full vector $\q$ can be
measured by various distances. A natural one is the
Euclidean distance, with square
\[
\bigl\| \te- \te' \bigr\|^2 = \sum_{i=1}^{n}
\bigl(\te_i-\te'_i\bigr)^2.
\]
If the indices of the $p_n$ nonzero coordinates of a
vector in the model $\ell_0[p_n]$ were known a priori, then the vector
could be estimated with mean square error of the order $p_n$. In
\cite{djhs92} it is shown that, as $n,p_n\ra\infty$ with $p_n=o(n)$,
\[
\inf_{\hat{\te}} \sup_{\te\in\ell_0[p_n]} P_{n,\te} \bigl\|\hat{\te}-\te
\bigr\|^2 = 2p_n \log(n/p_n) \bigl(1+o(1) \bigr).
\]
Here the infimum is taken over all estimators $\hat\q=\hat\q(X)$, and
$P_{n,\q}$ denotes taking the expectation under the assumption that
$X$ is $N_n(\q,I)$-distributed. In other words, the square minimax rate
over $\ell_0[p_n]$ is $p_n\log(n/p_n)$, meaning that the unknown
identity of the nonzero means needs to lead only to a logarithmic
loss.

The Bayesian approach is
presumably adopted for the intuition provided by prior modeling,
and is not necessarily directed at attaining minimax rates.
However, for theoretical investigation, it is natural to
take the minimax rate as a benchmark, and it is of
particular interest to know which priors yield a posterior
distribution that concentrates most of its mass on balls around $\q_0$
of square radius of order $p_n\log(p_n/n)$, or close
relatives as $p_n(\log n)^r$ that loose (only)
a logarithmic factor. A second main result of
the paper is that the minimax rate is attained for many
combinations of priors. It suffices that the priors $\pi_n$ decrease
exponentially with dimension, and give
sufficient weight to the true level of sparsity: for some $c>0$,
%
%e1.2 #&#
%
\begin{equation}
\label{EqLBOnPin} \pi_n(p_n)\gtrsim\exp
\bigl(-c p_n\log(n/p_n) \bigr).
\end{equation}
Furthermore, the priors on the nonzero coordinates should have
tails that are not lighter than Laplace, and satisfy a number of other
technical properties. If inequality~(\ref{EqLBOnPin}) fails,
then the rate of contraction may be slower than minimax; we show
that it is not slower than $\log(1/\pi_n(p_n) )$.
[The word ``contraction'' is in line
with other literature on nonparametric Bayesian procedures;
with the present choice of metrics (which grow with $n$)
the rates actually increase to infinity.]

More generally, we consider reconstruction relative to the
$\ell^q$ metric for \mbox{$0<q\le2$}, defined (without $q$th root) by
%
%e1.3 #&#
%
\begin{equation}
\label{EqDefdq} d_q\bigl(\te,\te'\bigr) = \sum
_{i=1}^{n} \bigl|\te_i-
\te'_i\bigr|^q.
\end{equation}
For $q<2$ this ``metric'' is more sensitive to
small variations in the coordinates than the square Euclidean metric,
which is $d_2$. (For $q\le1$ the definition gives a true metric $d_q$;
for $1< q\le2$ it does not.)
From~\cite{djhs92} the minimax rate over $\ell_0[p_n]$ for $d_q$
is known to be of the order
%
%e1.4 #&#
%
\begin{equation}
\label{EqDefrnq} r_{n,q}^*=p_n\log^{q/2}
(n/p_n).
\end{equation}
We show that the
posterior ``contraction'' rate attains this order under conditions as
in the preceding paragraph, and more generally characterize
the rate in terms of $\log(1/\pi_n(p_n))$.

Besides nearly black vectors, we consider rates of contraction
if $\q_0$ is in a weak $\ell_s$-ball. The minimax rate over
$m_s[p_n]$ relative
to $d_q$ is (see~\cite{DJ94})
%
%e1.5 #&#
%
\begin{equation}
\label{EqDefmunq} \mu_{n,s,q}^*=n \biggl(\frac{p_n}n
\biggr)^s\log^{(q-s)/2}(n/p_n).
\end{equation}
This is shown to be also the rate of posterior contraction under
slightly stronger conditions on the priors than before: the prior on
dimension must decrease slightly faster than exponential. Under the
same conditions we also show that the posterior distribution
has exponential concentration, and therefore
contracts also in the stronger sense of (any, Euclidean) moments.

A summary of these results is that good priors for the dimension
decrease at exponential or, perhaps better,
slightly faster rate, and good priors on
the nonzero means have tails that are heavier than Laplace. We also
show that priors with lighter tails, such as the Gaussian, attain
significantly lower contraction rates at true parameter vectors $\q_0$
that are not close to the origin.

The structure of the article is as follows. In Section
\ref{SectionMainResults} we state the main concentration results.
A practical algorithm, simulations and some pictures are presented in Section
\ref{SectionSimulations}. Proofs are gathered at
the end of the paper and in the supplementary Appendix~\cite{spa-appendix}.

%s1.1 #&#
\subsection{Notation}
We denote by $a\wedge b$ and $a\vee b$ the
minimum and maximum of two real numbers $a,b$, and
write $a\lesssim b$ if $a\le C b$ for a universal constant $C$.
The notation $\triangleq$ means ``equal by definition to.''
%We denote by $\Phi$ the
%cumulative distribution function of a standard normal variable and
%let $\bar{\Phi}=1-\Phi$.
We call \textit{support} of a vector $\te=(\q_1,\ldots,\q_n)\in\RR^n$
the set of indices of nonzero coordinates,
and denote this by $S_\te=\{i\in\{1,\ldots,n\}\dvtx\te_i\neq0\}$.
We set $\q_S=(\q_i\dvtx i\in S)$, and
let $|S|$ be the cardinality of a set $S\subset\{1,\ldots,n\}$.

%s2 #&#
\section{Main results}
\label{SectionMainResults}
Throughout the paper we consider a prior $\Pi_n$ on $\RR^n$
constructed in three steps:
\begin{longlist}[(P3)]
\item[(P1)] A \textit{dimension} $p$ is chosen according to a prior
probability measure $\pi_n$ on the set $\{0,1,2,\ldots, n\}$.
\item[(P2)] Given $p$ a subset $S\subset\{1,\ldots, n\}$ of size
$|S|=p$ is chosen uniformly at random from the ${n\choose p}$
subsets of size $p$.
\item[(P3)] Given $(p,S)$ a vector $\q_S=(\te_i\dvtx i\in S)$ is chosen
from a probability distribution with Lebesgue density $g_S$ on $\RR^p$ (if
$p\ge1$), and this is extended to $\te\in\RR^n$ by setting the
remaining coordinates $\q_{S^c}$ equal to 0.
\end{longlist}
For simplicity we use the same density $g_S$ for every set
of a given dimension $|S|$, and will denote this also by $g_{|S|}$.
We also assume that the prior on dimension is positive, that is
$\pi_n(p)>0$ for any integer $p$.

Given the prior $\Pi_n$, Bayes's rule yields the \textit{posterior
distribution} $B\mapsto\Pi_n(B\given X)$,
the conditional distribution of $\q$ given $X$ if the
conditional distribution of $X$ given $\q$ is taken equal to the
normal distribution $N_n(\q,I)$. The probability $\Pi_n(B\given X)$
of a Borel set
$B\subset\RR^n$ under the posterior distribution can be written
%
%e2.1 #&#
%
\begin{equation}
\label{EqPosterior} \frac{ {\sum_{p=0}^n{\pi_n(p)}/{{n\choose
p}}\!\sum_{|S|=p}
\int_{(\q_S,0)\in B}\!
\prod_{i\in S} \phi(X_i-\te_i)\!\prod_{i\notin S}\phi(X_i) g_S(\q_S) \,
d\te
_S}} {
{\sum_{p=0}^n{\pi_n(p)}/{{n\choose p}}\!\sum_{|S|=p}\int\!
\prod_{i\in S} \phi(X_i-\te_i)\!\prod_{i\notin S}\phi(X_i) g_S(\q_S) \,
d\te_S}}.\hspace*{-35pt}
\end{equation}
Here $(\q_S,0)$ is the vector in $\RR^n$ formed by adding coordinates
$\q_i=0$ to $\q_S=(\q_i\dvtx i\in S)$, at the positions left open
by $S\subset\{1,\ldots, n\}$ (in the correct order
of the coordinates and not at the end, as the notation suggests).
This expression is somewhat unwieldy; we consider computation in
Section~\ref{SectionSimulations}.
%The posterior distribution
%can be marginalized to obtain the posterior distribution of
%an arbitrary transformation $f(\q)$ by choosing the set $B$
%equal to $B=\{\q: f(\q)\in B_0\}$.

The posterior distribution is a random probability distribution on
$\RR^n$, which we study under the assumption that the vector
$X=(X_1,\ldots, X_n)$ is distributed\vadjust{\goodbreak} according to a multivariate
normal distribution with mean vector $\te_0$ and covariance matrix the
identity. We let $P_{n,\te_0} T$ denote the expected value of a
function $T=T(X)$ under this distribution.

We shall be interested in two aspects of the posterior distribution:
its dimensionality and its ability to recover the mean vector $\q$.
Because the conditions are simpler in the case that
the nonzero coordinates are independent under the prior,
in the first two results we assume that the densities $g_S$ in
(P3) are of product form. Concrete examples of priors as in (P1) and (P3)
that satisfy the conditions imposed in the theorems
are given in Section~\ref{SectionExamples}.

%s2.1 #&#
\subsection{Dimensionality}
In the context of $\ell_0[p_n]$-classes, we say that the prior $\pi_n$
on dimension has \textit{exponential decrease} if, for some constants $C>
0$ and $D<1$,
%for some constants $C>0$ and $D<1$,
%
%e2.2 #&#
%
\begin{equation}
\label{priordim} \pi_n(p) \le D \pi_n(p-1),\qquad
p > Cp_n.
\end{equation}
If the condition is also satisfied with $C=0$, we say that the prior on
dimension has \textit{strict} exponential decrease.
%For weak classes, the definition is analogous,
%but $p_n$ should be replaced by $p_n^*=n(p_n/n)^s\log^{-s/2}(n/p_n)$.
%Our condition is in fact slightly weaker, see Section
%
%th2.1 #&#
%
\begin{theorem}[(Dimension)]
\label{TheoremDimension}
If $\pi_n$ has exponential decrease~(\ref{priordim})
and $g_S$ is a product of $|S|$ copies of a univariate
density $g$,
with mean zero and finite second moment,
then there exists $M>0$ such that, as $p_n, n\ra\infty$,
\[
\sup_{\te_0\in\ell_0[p_n]} P_{n,\te_0} \Pi_n \bigl(\te\dvtx
|S_\te|>Mp_n \given X \bigr) \to0.
\]
\end{theorem}

%It is to be expected
For reasonable priors, we may hope that the posterior distribution
spreads mass in
the $p_n$-dimensional subspace that supports a true mean vector
$\q_0\in\ell_0[p_n]$. The theorem shows that the posterior
distribution ``overshoots'' this space by subspaces of dimension at
most a multiple of $p_n$. Because the overshoot can have a random
direction, this does not mean that the posterior distribution
concentrates overall on a fixed $Mp_n$-dimensional subspace. The
theorem shows that it concentrates along $Mp_n$-dimensional coordinate
planes, but its support will be far from convex.

Obviously the posterior distribution will concentrate on low-dimensional
subspaces if the higher-dimensional spaces receive
little mass under the prior $\pi_n$. By the theorem,
exponential decrease is sufficient. The next
step is to show that exponential decrease is not too harsh:
it is compatible with good reconstruction of the full mean vector $\q$.
This then, of course, requires a lower bound on the prior mass
given to the spaces of ``correct'' dimension; for instance, see
(\ref{EqLBOnPin}).

%s2.2 #&#
\subsection{Recovery}
Good recovery requires also appropriate prior densities $g_S$ on the
nonzero coordinates. Because the statistical problem of recovering
$\te$ from a $N_p(\q,I)$ distributed observation is equivariant in
$\te$, we may hope that the location of the nonzero coordinates of
$\te_0$ does not play a role in its recovery rate. The non-Bayesian
procedures considered in, for instance,~\cite{Golubev02} indeed fulfill
this expectation. However, a Bayesian procedure (with proper priors)
necessarily favors certain regions of the parameter space. Depending
on the choice of priors $g_S$ in~(P3), this may lead to a shrinkage
effect, even in the ``average'' recovery of the parameter as
$n\ra\infty$, yielding suboptimal behavior for true parameters
$\te_0$ that are far from the origin. This shrinkage effect can be
prevented by choosing priors $g_S$ with sufficiently heavy tails.
% Remark: in the work by Jiang and Zhang(2009), the ``distribution of
%the
% theta_i's" is estimated via a NPMLE which is equivariant. Their final
%estimator
% (which is empirical Bayes rather than Bayes) is then also equivariant.

Again we first consider the case of independent coordinates.
In the following theorem
we assume that $g_S$ is a product of $|S|$ densities of the form
$e^h$, for a function $h\dvtx\RR\to\RR$ satisfying
%
%e2.3 #&#
%
\begin{equation}
\label{EqProductPrior} \bigl|h(x)-h(y) \bigr|\lesssim1+|x-y|\qquad \forall x,y\in\RR.
\end{equation}
This covers all densities $e^h$ with a uniformly Lipshitz function
$h$, such as the Laplace and Student densities. (For the Student
density the following theorem assumes more than 2 degrees
of freedom to ensure also
finiteness of the second moment.) It also covers other smooth
densities with polynomial tails, and densities of the form
$c_\alpha e^{-|x|^\al}$ for some $\al\in(0,1]$, which have a
function $h$ that
is bounded in a neighborhood of the origin and uniformly Lipschitz
outside the neighborhood. On the other hand the standard normal
density is ruled out. In Theorem~\ref{TheoremLB} we shall see that
this indeed causes a shrinkage effect.

Recall definition~(\ref{EqDefdq}) of the (square) distance $d_q$.
%
%th2.2 #&#
%
\begin{theorem}[(Recovery)]
\label{TheoremConcentrationIndependent}
If $\pi_n$ has exponential decrease~(\ref{priordim})
and $g_S$ is a product of $|S|$ univariate
densities of the form $e^h$ with mean zero and finite second moment
and $h$ satisfying~(\ref{EqProductPrior}),
then for any $q\in(0,2]$,
for $r_n$ satisfying
%
%e2.4 #&#
%
\begin{equation}
\label{EqConditionRate} r_n^2 \geq\bigl\{
p_n\log(n/p_n) \bigr\} \vee\log\frac1{
\pi_n(p_n)}
\end{equation}
and sufficiently large $M$,
as $p_n,n\to\infty$ such that $p_n/n\ra0$,
\[
\sup_{\te_0\in\ell_0[p_n]} P_{n,\te_0} \Pi_n \bigl(\te\dvtx
d_q(\q,\q_0)> M r_n^q
p_n^{1-q/2} \given X \bigr) \to0.
\]
\end{theorem}

For $q=2$ the theorem refers to the square Euclidean distance $d_2$,
and asserts that the posterior distribution contracts at the rate
$r_n^2$, uniformly over $\ell_0[p_n]$. The first inequality in
(\ref{EqConditionRate}) says that this rate is (of course) not faster
than the minimax rate $r_{n,2}^*=p_n\log(n/p_n)$. The second shows
that it is also limited by the amount of prior mass $\pi_n(p_n)$ put on
the true dimension. If this satisfies~(\ref{EqLBOnPin}), then
$\log(1/\pi_n(p_n))\lesssim r_{n,2}^*$ and the
rate $r_n^2$ is equal to the minimax rate.

Condition~(\ref{EqLBOnPin}) for every $p_n$ leaves a free margin of a
$\log(n/p_n)$-term over just exponential decrease of the prior
$\pi_n$. If the decrease is still faster than~(\ref{EqLBOnPin}),
then the rate of contraction may
be slower. For instance, for $\pi_n(p)\asymp\exp(-p^\alpha)$, for some
$\alpha>1$, the rate for the square Euclidean distance given by the
theorem is not better than $p_n^\alpha$, which is much bigger than
$r_{n,2}^*$.
In contrast, for $\alpha=1$ the theorem gives the minimax rate.
%Is this sharp?

For $q\in(0,2)$ we can make similar remarks. The minimax rate
$r_{n,q}^*$ over $\ell_0[p_n]$ for $d_q$ is given in~(\ref{EqDefrnq}).
Because
\[
\bigl(r_{n,2}^*\bigr)^{q/2}p_n^{1-q/2}=r_{n,q}^*,
\]
the theorem shows
contraction of the posterior distribution relative to $d_q$ at the
minimax rate $r_{n,q}^*$ over $\ell_0[p_n]$ under the\vspace*{1pt} same conditions
that it gives the minimax rate $r_{n.2}^*$ for $d_2$:
(\ref{EqLBOnPin}) suffices. Furthermore, if there is less prior mass
at $p_n$, then the rate of contraction will be slower.

In the case that $0<q<1$ the result is surprising at first when
compared to the finding in~\cite{JohnSilv04} that the posterior
\textit{median}, or more generally so-called ``strict-thresholding
rules,'' attain the convergence rate $r_{n,q}^*$, but the posterior
\textit{mean} converges at a \textit{strictly slower} rate (even when
$\q_0=0$; see Section 10 in~\cite{JohnSilv04} and the remark below). By
the preceding theorem the \textit{full} posterior distribution
\textit{does} contract at the optimal rate $r_{n,q}^*$, for any
$0<q<2$. This is true in particular for the case of binomial priors on
dimension considered in~\cite{JohnSilv04} with the ``best possible''
(oracle) choice $\al_n=p_n/n$.

The slower convergence of the posterior mean relative to the
contraction of the full posterior distribution is made possible by the
fact that $d_q$-balls have astroid-type shapes for $0<q<1$, and differ
significantly from their convex hull if $n$ is large. The posterior
mean, which is in the convex hull of the support of the posterior, can
therefore be significantly farther in $d_q$-distance from $\q_0$ than
the bulk of the distribution. By Theorem~\ref{TheoremDimension} only
few coordinates outside the support of $\q_0$ are given nonzero values
by the posterior. However, the corresponding indices are random and
\textit{on average} spread over $\{1,2,\ldots,n\}$, which makes that
the posterior mean at a fixed coordinate is typically nonzero. Adding
up all small errors in $\ell^q$ typically gives a much higher total sum
for $q<1$ than for $q\ge1$. In contrast the posterior median does not
suffer from this averaging effect.

The posterior measure thus provides a unifying point of view on the
considered objects. In this perspective for $0<q<1$ the posterior mean
is a bad representation of the full posterior measure.
%
%re2.3 #&#
%
\begin{remark}
From the arguments exposed in~\cite{JohnSilv04}, it is not hard to
check that the posterior mean generally fails to attain the minimax
rate over $\ell_0[p_n]$ relative to $d_q$ for $0<q<1$. Let us consider
the case of $\ell_0[p_n]$ classes with $\te_0=0$. With the notation of
\cite{JohnSilv04}, the posterior mean $\tilde{\mu}(x,\al_n)$ with data
$X_1=x$ for the binomial prior on dimension with parameters
$(n,\al_n)$ satisfies $|\tilde{\mu}(x,\al_n)| \ge C|x|\al_n$, by the
same reasoning as\vadjust{\goodbreak} in the last display of page 1647 in
\cite{JohnSilv04} (the weight parameter $\hat{w}$ is fixed here and
equals $\al_n$). Hence the $\ell^q$-power loss $\sum_i
P_{n,\te_0}|\te_{0,i}-\tilde{\mu}(X_i,\al_n)|^q$ when $\te_0=0$ is
bounded from below by a constant times $n\al_n^q$. Thus, even for the
``oracle'' parameter $\al_n=p_n/n$, this is much above the minimax
risk for any $0<q<1$.
\end{remark}

%s2.3 #&#
\subsection{Dependent priors}
The preceding theorems are also true for priors that render the
coordinates $\q_i$ dependent. In the remaining theorems
we assume that the densities $g_S$ in
(P3) satisfy the conditions, for every $S'\subset S\subset\{1,\ldots,n\}$
and a universal constant $c_1$,
%
%e2.5 #&#
%e2.6 #&#
%
\begin{eqnarray}
\label{EqGOne}
\log g_S(\q)-\log g_S\bigl(\q'\bigr)&\le&
c_1|S|+\tfrac1{64} {\bigl\|\q-\q'\bigr\|^2}\qquad \forall
\q,\q'\in\RR^S,
\\
\label{EqGTwo}
\bigl|\log g_S(\q)-\log g_{S'}(\pi_{S'}\q) \bigr| &\le&
c_1|S| +\tfrac1{64} {\|\pi_{S-S'}\q\|^2}\qquad
\forall\q\in\RR^S.
\end{eqnarray}
Here $\pi_S\dvtx\RR^n\to\RR^S$ is the projection defined by
$\pi_S\q=\q_S=(\q_i\dvtx i\in S)$.
(The constant 64 corresponds to the constant 32 in Lemma~\ref{LemmaTestPhi},
but has no special significance and can be improved.)

For a partition $S=S_1\cup S_2$, we denote by $\q=(\q_1,\q_2)$ the
corresponding partition of $\q\in\RR^S$ and
by $g_{S_1,S_2}(\q_1,\q_2)=g_S(\te)$ the corresponding density.
%$\q_2\mapsto g_{S_2|S_1}(\q_2\given\te_1)\propto g_S(\q_1,\q_2)$
%the conditional density of $\te_2$ given $\te_1$ resulting from $g_S$.
In the next theorem we assume that there exist $C, m_1>0$
%these conditional densities are centered,
%and have uniformly bounded marginal second moments: for some
%constant $m_2$ and every $S_1\subset S_2^c$,
and, for any $S_2$, probability densities
$\ga_{S_2}$ on $\RR^{S_2}$, such that for any $\te_2\in\RR^{S_2}$
and $S_1\subset S_2^c$,
%
%e2.7 #&#
%
\begin{equation}
\label{EqBoundConditionalDensity} \sup_{\te_1\in\RR^{S_1}} \frac{
g_{S_1,S_2}(\te_1,\te_2) }{ g_{S_1}(\te_1) }\le
Cm_1^{|S_1|+|S_2|}\ga_{S_2}(\te_2).
\end{equation}
This condition expresses that the ``mixing
between the coordinates within a given subspace'' is not too important.

Examples are given in Section~\ref{SectionExamples}.
%
%th2.4 #&#
%
\begin{theorem}[(Recovery)]
\label{TheoremConcentrationDependent}
Suppose $\pi_n$ has strict exponential decrease, that is, satisfies
(\ref{priordim}) with $C=0$ and some $D>0$.
The assertions of Theorems~\ref{TheoremDimension}
and~\ref{TheoremConcentrationIndependent} are also true if
the densities $g_S$ are not product densities, but general
densities with finite second moments that satisfy
(\ref{EqGOne}),~(\ref{EqGTwo}) and~(\ref{EqBoundConditionalDensity})
with \mbox{$Dm_1<1$},
%for $D$ the constant in \eqref{priordim}
and $m_1$ the constant in
(\ref{EqBoundConditionalDensity}).
\end{theorem}

% \Pi(\cdot\given\bar\theta_1) \Pi_n(S_2\given S_1) \pi_n(k)

%s2.4 #&#
\subsection{Complexity priors}
The next results are designed for application to the particular
priors $\pi_n$ of the form, for positive constants $a,b$,
%
%e2.8 #&#
%
\begin{equation}
\label{priorone} \pi_n(p) \propto e^{-ap\log(bn/p)},
\end{equation}
where $\propto$ stands for ``proportional to.'' Because $e^{p\log
(n/p)}\le{n \choose p} \le e^{p\log(ne/p)}$, this prior
is inversely proportional to the number of models of size $p$,
a quantity that could be viewed as the \textit{model complexity}
for a given dimension $p$. Thus this prior
appears particularly suited to the purpose of ``downweighting the
complexity.'' Forgetting about the extra component $g_S$ of
the prior, we can also consider it an analog of the penalty
``$2p\log(n/p)$'' used in model selection in this context by (e.g.)
Birg\'e and Massart in~\cite{BirgeMas01}. Every particular
model with support $S$ of size $|S|=p$ receives prior probability
bounded below and above by expressions of the type
$e^{-a_1 p\log(b_1n/p)}$ from this prior.

Because the \textit{complexity prior}~(\ref{priorone}) has exponential
decrease~(\ref{priordim}) when $b > 1+e$ and satisfies~(\ref{EqLBOnPin}),
Theorems~\ref{TheoremDimension} and~\ref{TheoremConcentrationDependent}
(or Theorem~\ref{TheoremConcentrationIndependent})
show that the corresponding posterior distribution
concentrates on low-dimensional spaces and attains the
minimax rate of contraction over $\ell_0[p_n]$ relative to (any) $d_q$
if combined with densities $g_S$ satisfying the conditions of
Theorem~\ref{TheoremConcentrationDependent}. The following theorem
relaxes the condition on $g_S$ and gives a more precise result
on the contraction of the posterior measure.

The theorem applies more generally to priors on dimension
satisfying the upper bound, for some $a,b>0$, and
every $p\in\{0,1,\ldots, n\}$,
%
%e2.9 #&#
%
\begin{equation}
\label{prioroneonesided} \pi_n(p)\lesssim e^{-ap\log(bn/p)}.
\end{equation}

%th2.5 #&#
%
\begin{theorem}[(Recovery)]
\label{TheoremMainExponentialDecrease}
If the densities $g_S$ have finite second moments, satisfy (\ref
{EqGOne}) and~(\ref{EqGTwo}) for some constant $c_1$,
and the priors $\pi_n$ satisfy~(\ref{prioroneonesided})
for some $a\ge1$ and $b\ge e^{7+2c_1}$,
then, for $r_n$ satisfying~(\ref{EqConditionRate}),
for any $1\le p_n\le n$ and any $r\ge1$,
\[
\sup_{\te_0\in\ell_0[p_n]}P_{n,\q_0}\Pi_n \bigl(\q\dvtx\|\q-
\q_0\| >45r_n+10r\given X \bigr) \lesssim e^{- r^2/10}.
\]
\end{theorem}

Consistent with the preceding findings, the posterior distribution concentrates
on Euclidean balls of radius of the order $r_n$ around $\q_0$. In addition
the theorem shows that its ``tail'' is sub-Gaussian,
uniformly in $n$ and uniformly over $\ell_0[p_n]$.
As one consequence, for every $l\in\NN$,
\[
P_{n,\q_0}\int\|\q-\q_0\|^l \,d
\Pi_n(\q\given X) \lesssim r_n^l.
\]
By Jensen's inequality, this in turn implies the following corollary.
%$$P_{n,\q_0}\Bigl\|\int\q d\Pi_n(\q\given X)-\q_0\Bigr\|^l
%Thus the risk of the posterior mean $\int\q d\Pi_n(\q\given X)$
%as a point estimator of $\q_0$
%is of the order $r_n$, relative to every polynomial loss function.
%
%co2.1 #&#
%
\begin{corollary}[(Posterior mean)]
Under the conditions of Theorem~\ref{TheoremMainExponentialDecrease},
%the posterior mean $\hat{\q}^{PM}_n= \int\q d\Pi_n(\q\given X)$
%satisfies,
%
\[
\forall l\in\mathbb{N}\qquad \sup_{\te_0\in\ell_0[p_n]} P_{n,\q_0} \biggl\| \int\q \,d
\Pi_n(\q\given X) -\q_0 \biggr\|^l \lesssim
r_n^l.
\]
\end{corollary}
The posterior mean $\int\q\,d\Pi_n(\q\given X)$
as a point estimator of $\q_0$ has a risk of the order $r_n$, relative
to every polynomial loss function. In particular, it is rate-minimax
over $\ell_0[p_n]$ for the squared $\ell_2$-risk.

%Other functionals of the posterior measure can also be of interest.
%For instance,
The posterior coordinate-wise median considered in the simulation study
below is another interesting functional of the posterior measure. Under
the conditions of Theorem~\ref{TheoremMainExponentialDecrease} and
(\ref{priorone}), the posterior coordinate-wise median is rate-minimax
over $\ell_0[p_n]$, \textit{for any}
$d_q$-distance, $0<q\le2$; see~\cite{spa-appendix}.\vadjust{\goodbreak}

The theorem, with its explicit bound, is also the basis for
results on the concentration of the posterior distribution when
the true vector is in a weak $m_s[p_n]$-class. Results for the posterior
mean and $\ell_2$-risk can be obtained as above as a consequence.
%
%th2.6 #&#
%
\begin{theorem}[(Recovery, weak class)]
\label{TheoremWeakClass}
If the densities
$g_S$ have finite second moments, satisfy
(\ref{EqGOne}) and~(\ref{EqGTwo}) for some constant $c_1$,
and the priors $\pi_n$ satisfy~(\ref{prioroneonesided})
for some $a\ge1$ and $b\ge e^{7+2c_1}$, then, for $r_n$ satisfying
\[
r_n^2 =\min_{1\le p\le n} \biggl[ \frac{sn^{2/s}}{2-s}
\biggl(\frac1p \biggr)^{2/s-1} \biggl(\frac
{p_n}n
\biggr)^2 \vee p\log\frac np \vee\log\frac1{\pi_n(p)}
\biggr]
\]
for any $1\le p_n\le n$, $s\in(0,2)$ and any $r\ge1$,
\[
\sup_{\te_0\in m_s[p_n]}P_{n,\q_0}\Pi_n \bigl(\q\dvtx\|\q-
\q_0\|> 80r_n+20r\given X \bigr) \lesssim e^{- r^2/10}.
\]
\end{theorem}

For the ``complexity prior'' $\pi_n$ given by~(\ref{priorone})
the third term $\log(1/\pi_n(p) )$ in the minimum
defining it is smaller than a multiple of the second term, and hence
can be omitted. The minimum can then be determined by equating
the first two terms, leading to
% n^{2/s}\left(\frac{p_n}{n}\right)^2 \gtrsim\log n
%This amounts to solve in $K_n$ the equation
%$$ n^{\frac{2}{s}}\left(\frac{p_n}{n}\right)^2 K_n^{1-\frac{2}{s}} =
%K_n\log(n/K_n).$$
%
%e2.10 #&#
%
\begin{equation}
\label{EqDefpnstar} p_n^*\asymp n(p_n/n)^s/
\log^{s/2}(n/p_n).
\end{equation}
If $p_n^*\gtrsim1$, then this value is eligible in the minimum,
and the first and second terms evaluated at $p_n^*$ are of the same order,
given by
\[
r_n^2\asymp n \biggl(\frac{p_n}n
\biggr)^s\log^{1-s/2}\frac n{p_n}.
\]
This in fact is the minimax rate $\mu_{n,s,2}^*$ for the square Euclidean
metric $d_2$ over the class $m_s[p_n]$; see~(\ref{EqDefmunq}).
Thus the complexity priors combined with densities $g_S$
satisfying~(\ref{EqGOne}) and~(\ref{EqGTwo}) [in particular,
product densities satisfying~(\ref{EqProductPrior})]
yield contraction at the minimax rate
over both the nearly black vectors $\ell_0[p_n]$ and the weak $m_s[p_n]$
classes. For priors on dimension that are significantly smaller
than the complexity priors, the third term in the minimum
must be taken into account, and the rate of contraction is smaller than minimax.

The condition $p_n^*\gtrsim1$ is satisfied as soon as
the sparsity coefficient $p_n/n$ is not too small. If the signal is
very sparse
and has $p_n^*\ll1$, then the minimum in the definition of
$r_n^2$ is taken at $p\sim1$, leading to a squared rate of the order
$\log{n}$.
This is within a constant of the rate
achieved by hard thresholding in that case.
%One should have results for the euclidian distance with an extra log
%for
%all priors for which the extension of Thm 1 (spa6) works. Since the
%convergence is exponentially fast, we should also have the behavior of
%the mean for those classes.

The previous result extends under slightly stronger conditions to
$d_q$-distances with $q>s$. Furthermore, the following theorem
shows that $p_n^*$ is indeed an upper bound on the dimensionality\vadjust{\goodbreak}
of the posterior distribution.
For simplicity we only state the result in the case
of complexity priors. Recall that $\m_{n,s,q}^*$, given
in~(\ref{EqDefmunq}), denotes the minimax rate
over the class $m_s[p_n]$ relative to $d_q$.
%
%th2.7 #&#
%
\begin{theorem}[(Dimensionality, recovery, weak class)] \label{thmpq}
Suppose the densities $g_S$ have finite second moments, satisfy (\ref
{EqGOne}),~(\ref{EqGTwo}) and~(\ref{EqBoundConditionalDensity}),
and $\pi_n$ satisfies~(\ref{priorone}) for sufficiently
large $a\ge1$ and $b>e$. Then
for any $s\in(0,2)$, any $q\in(s,2)$ and any $(p_n)$ such that
$p_n/n\to0$ and $p_n^*$ given by~(\ref{EqDefpnstar}) is bounded
away from 0, for a sufficiently large constant $M$,
\begin{eqnarray*}
\sup_{\te_0\in m_s[p_n]}P_{n,\q_0} \Pi_n \bigl(\q\dvtx
|S_\q|> Mp_n^* \given X \bigr)&\to&0,
\\
\sup_{\te_0\in m_s[p_n]}P_{n,\q_0} \Pi_n \bigl(\q\dvtx
d_q(\q,\q_0)> M\mu^*_{n,s,q} \given X \bigr)&\to&0.
\end{eqnarray*}
\end{theorem}

%s2.5 #&#
\subsection{Examples}
\label{SectionExamples}
In this section we discuss examples of priors on dimension
$\pi_n$ and prior densities $g_S$ on the nonzero coordinates that
satisfy the conditions of the preceding theorems.
%
%ex2.1 #&#
%
\begin{example}[(Independent Dirac mixtures)]
Consider the prior on $\q=(\q_1,\ldots, \q_n)\in\RR^n$
corresponding to sampling the coordinates $\te_i$ independently from a mixture
$(1-\al)\delta_0 + \al g$ of a Dirac measure at 0 and a univariate
density~$g$,
for a given $\al\in(0,1)$. The coordinates of $\q$ are then independently
zero with probability $1-\al$, and hence the dimension of the
model is binomially distributed with parameters $n$ and $\al$.
Furthermore, the nonzero coordinates are distributed according to the
product of copies of $g$. Thus this prior fits in our set-up, with
$\pi_n$ the binomial$(n,\al)$-distribution and $g_S$ a product
density.

For a fixed $\al$ the coordinates $\q_i$ are independent,
under both the prior and the posterior distribution. Furthermore,
the posterior distribution of $\q_i$ depends on $X_i$ only.

This prior is considered in~\cite{GeorgeFoster}
and~\cite{JohnSilv04}, in combination with a
Gaussian or a heavy tailed density $g$, respectively.
In the next section we show that Gaussian priors are
deficient if the nonzero coordinates of the signal are
large. The authors of~\cite{JohnSilv04} propose to use the
coordinatewise posterior median (or another univariate point
estimator) for estimating $\q$, with the weight parameter $\alpha$
set by
a thresholded empirical Bayes method: the parameter is chosen equal to
the maximum likelihood estimator of $\al$ based on the marginal
distribution of $X$ in the Bayesian set-up (i.e., with $\q$ integrated
out but with fixed $\al$) subject to the constraint that the resulting
posterior median (after plugging in $\hat\al$) given an observation in
the interval $[-\sqrt{2\log n},\sqrt{2\log n}]$ is zero. The authors
show that the resulting point estimator works remarkably well, in a
minimax sense, for various metrics and sparsity classes.\vadjust{\goodbreak}

A natural Bayesian approach is to put a prior on $\al$, which yields a
mixture of binomials as a prior $\pi_n$ on the dimension of the
model. The independence of the coordinates $\q_i$ is then lost.
We discuss this prior further in the following example.
\end{example}
%
%ex2.2 #&#
%
\begin{example}[(Binomial and beta-binomial priors)]
\label{priorhie}
The binomial $(n,\alpha_n)$ distribution as the prior
$\pi_n$ on dimension gives an expected dimension of $n\alpha_n$.
%$$ \pi^{B,\al}_{n}(k) = {n \choose k} \al^k (1-\al)^{n-k}.$$
In the sparse setting a small value of $\alpha_n$ is therefore natural.
If the sparsity parameter $p_n$ were known, we could
consider the choice $\al_n=p_n/n$; we shall refer to the corresponding law
as \textit{oracle binomial prior}.
%(it can be seen as a Bayes analog of the oracle thresholding $\sqrt{2

Assume that $p_n\ra\infty$ with $p_n/n\ra0$.
The binomial prior has exponential decrease~(\ref{priordim})
if $\al_n \lesssim p_n/n$. The oracle binomial prior
$\alpha_n\asymp p_n/n$ is at the upper end of this range, and also
satisfies~(\ref{EqLBOnPin}), and thus yields the minimax
rate of contraction.
The choice $\alpha_n=1/n$ yields
$\log\pi_n(p_n)$ of the order $-p_n\log p_n$, and hence attains
the minimax rate if $p_n$ is of the order $n^{a}$, $a<1$; for larger
$p_n$ it may miss the minimax rate by a logarithmic factor.

A natural Bayesian strategy is to view the unknown ``sparsity''
parameter $\alpha$ as a hyperparameter and put a prior on it. The classical
choice is the Beta prior, leading to the hierachical scheme
$\al\sim\operatorname{Beta}(\k,\l)$ and $p\given\al\sim
\operatorname{binomial}(n,\al)$,
which corresponds to the following prior on $p$:
\[
\pi_n(p)=\pmatrix{n\cr p}\frac{B(\k+p,\l+n-p)}{B(\k,\l)} \propto\frac
{\Gamma(\k+p)\Gamma(\l+n-p)}{p!(n-p)!}.
\]
The mean dimension is $n\k/(\k+\l)$, which suggests to choose the
hyper parameters of the Beta distribution so that $\k/(\k+\l)$ is in
the range $(c/n, C p_n/n)$. It is easy to verify that
the prior has exponential decrease~(\ref{priordim}), with $C=1$,
if $(\k-1)/p_n<D (\l-1)/(n-p_n+1)+D-1$. This suggests to choose
small $\k$ and large $\l$, thus giving a small variance to the
Beta distribution.
%$$\frac{\pi_n(p)}{\pi_n(p-1)}
%=\frac{\k+p-1}p\frac{n-p+1}{\l+n-p}=\frac{1+)\k-1)/p}{1+(

For $\k=1$ and $\l=n+1$ we obtain $\pi_n(p)\propto{2n-p\choose n}$.
Then $\pi_n(p)/\pi_n(p-1) = (n-p+1)/(2n-p+1)$, showing (strict) exponential
decrease~(\ref{priordim}), with $D=1/2$. By a binomial identity the norming
constant is equal to ${2n+1\choose n}$, so
\[
\pi_n(p)=\frac{(2n-p)(2n-p-1)\cdots(2n-p-n+1)} {
(2n+1)2n\cdots(2n+1-n+1)}\ge\biggl(1-\frac{p+1}{n+2}
\biggr)^n.
\]
For $p_n/n\ra0$, this gives $\pi_n(p_n)\gtrsim e^{-p_n(1+o(1))}$,
and hence~(\ref{EqLBOnPin}) is satisfied.
More generally, we may choose $\kappa=1$, $\la=\kappa_1n+1$,
which leads to $\pi_n(p)\propto{(\kappa_1+1)n-p \choose\kappa_1n}$.
The priors given by $\pi_n(p)\propto{{2n-p} \choose n}^{\kappa_1}$,
for some $\kappa_1>0$ are a further alternative.
\end{example}
%
%ex2.3 #&#
%
\begin{example}[(Poisson priors and hierarchies)]
The Poisson$(\alpha)$ distribution truncated to $\{0,1,\ldots,n\}$,
yields priors satisfying
\[
\pi_n(p)\propto\frac{e^{-\alpha}\alpha^p}{p!} \asymp C e^{-p\log
(p/\alpha)}e^p
\frac1{\sqrt p}\vadjust{\goodbreak}
\]
for $p\ra\infty$, by Stirling's approximation.
The mean is approximately $\alpha$, suggesting $\alpha$ in the range
$(1,cp_n)$. As $\pi_n(p)/\pi_n(p-1)=\alpha/p$, the prior
has exponential decrease~(\ref{priordim}) for $p\ge\alpha/D$.

If we put an exponential $(\l)$ hyperprior on $\alpha$, then
$\pi_n$ transforms into a shifted geometric distribution
(shifted $-1$ to have support starting at 0) with success
probability $\l/(1+\l)$. A Gamma hyperprior yields a
shifted negative binomial. For fixed hyper--hyper parameters
both are of the form $e^{-Cp}$ for some
constant $C$, and hence have exponential decrease, and satisfy~(\ref{EqLBOnPin}).
\end{example}
%
%ex2.4 #&#
%
\begin{example}[(Complexity prior)]
The prior $\pi_n(p)\propto e^{-ap\log(bn/p)}$ has exponential
decrease~(\ref{priordim}) for $b>1+e$ and satisfies~(\ref{EqLBOnPin}).
Theorems~\ref{TheoremMainExponentialDecrease},~\ref{TheoremWeakClass}
and~\ref{thmpq} show that this prior also gives sparsity and minimax recovery
of the parameter over weak $\ell_s$-classes. Although our results
do not show the opposite assertion that mere exponential decrease
is not enough for minimaxity on weak classes (while together
with~(\ref{EqLBOnPin}) it is enough for minimaxity over $\ell_0[p_n]$),
this might be a potential advantage of complexity priors
over the binomial and Poisson-based priors discussed
previously.
\end{example}
%
%ex2.5 #&#
%
\begin{example}[(Product prior)]
\label{exproduct}
Densities $g_S$ that are products of $|S|$ copies
of a univariate density with finite second moment
of the form $g=e^h$ for $h\dvtx\RR\to\RR$
a function that satisfies~(\ref{EqProductPrior}),
satisfy~(\ref{EqGOne}),~(\ref{EqGTwo}) and
(\ref{EqBoundConditionalDensity}).
In this sense Theorem~\ref{TheoremConcentrationDependent} is
a generalization of Theorem~\ref{TheoremConcentrationIndependent}.

To see this note that for a product density the function $g_S$
takes the form $g_S(\q)=\exp\{\sum_{i\in S}h(\q_i)\}$. Hence if
(\ref{EqProductPrior}) holds with proportionality constant~1, then
the left-hand side of~(\ref{EqGOne}) is bounded in absolute value by
\[
\sum_{i\in S}h(\q_i)-h\bigl(
\q_i'\bigr)\le|S|+\bigl\|\q-\q'
\bigr\|_1 \le|S|+\sqrt{|S|} \bigl\|\q-\q'\bigr\| \le5|S|+\frac1{64}\bigl\|
\q-\q'\bigr\|^2.
\]
Furthermore, the left-hand side of~(\ref{EqGTwo}) is bounded by
\[
\sum_{i\in S-S'} \bigl|h(\q_i) \bigr|
\le\bigl|S-S'\bigr| \bigl|h(0) \bigr| +\sum_{i\in S-S'}\bigl(1+|
\q_i|\bigr)\lesssim\bigl|S-S'\bigr|+\sum
_{i \in S-S'}|\q_i|.
\]
The $L_1$-norm of $(\q_i\dvtx i\in S-S')$ can be bounded by a linear
combination of $|S-S'|$ and the square $L_2$-norm, as before,
and hence the whole expression is
bounded by $C|S|+\|\pi_{S-S'}\q\|^2/64$, for some constant $C$.

Because a product density $g_S$
is a product of the marginal densities, the validity
of condition~(\ref{EqBoundConditionalDensity}) is clear.
\end{example}
%
%ex2.6 #&#
%
\begin{example}[(Weakly mixing priors)]
\label{wmprior}
For $h\dvtx\RR\to\RR$ a function satisfying~(\ref{EqProductPrior})
so that $e^h$ is integrable and $G\dvtx[0,\infty)\to\RR$ a
Lipschitz function that is bounded below,
consider, for $\te=(\te_1,\ldots,\te_p)$,
\[
g_p(\q) = a_p e^{\sum_{i=1}^p h(\te_i)-G(\|\te\|)},
\]
where $a_p$ is the normalizing constant.
An example is the prior, for $a>0$,
\[
g_p(\q) \propto\frac{e^{-\|\te\|_1}}{1+a^2\|\te\|^2}.
\]
In the Appendix~\cite{spa-appendix} it is shown that priors of this
form satisfy~(\ref{EqGOne}) and~(\ref{EqGTwo}).
%
% (###) Remark: If one takes $h(x)=-|x|$ and $\log G$ is $a$-Lipschitz
%with $a<1$,
% then "bounded away from zero" is not needed.
%
% (note that since $G$ is bounded away from zero, $|\log{a_p}|\leq Cp$
%, for some
% constant $C$).
Furthermore, it is shown that~(\ref{EqBoundConditionalDensity})
is also satisfied, with $m_1=(1+a)/(1-a)$ if $-h$ is the absolute value
of the identity function and the Lipschitz constant $a$ of $G$ is
strictly smaller than 1
[i.e., $ |G(s)-G(t) |\le a|s-t|$ for $a<1$].

%Because $g_p(\q)$ depends on the absolute values of
%the coordinates $\q_i$ only, the (conditional) densities
%are symmetric about 0 and hence centered.

Thus any prior of this form combined with any prior on dimension
that decreases exponentially such that $Dm_1=D(1+a)/(1-a)<1$,
for $D$ the constant in~(\ref{priordim}), gives
recovery at the minimax rate over $\ell_0[p_n]$,
by Theorem~\ref{TheoremConcentrationDependent}, and also
over $\ell_s[p_n]$ if combined with a complexity
prior on dimension satisfying the conditions
of Theorems~\ref{TheoremWeakClass} and~\ref{thmpq}.
For instance, the hierarchical binomial prior
$\pi_n(p)\propto{2n-p\choose n}$ in Example~\ref{priorhie}
has $D=1/2$ and hence $a<1/3$ suffices for contraction over
$\ell_0[p_n]$.
\end{example}
%
%[Laplace mixtures]\rm
%Choosing a scale parameter $\l$ from a
%Gamma $(\a,\b)$-distribution and given $\l$ generating
%the nonzero coordinates $\q_i$ from a Laplace distribution with
%scale $\l$, yields the densities, with $G$ the Gamma distribution,
%$$g_p(\q)=\int(\l/2)^p e^{-\l\|\q\|_1} dG(\l)=
%These densities have (conditional) mean zero,
%satisfy \eqref{EqGOne}-\eqref{EqGTwo}, and also

%s2.6 #&#
\subsection{Lower bounds} \label{seclb}
Condition~(\ref{EqProductPrior}) [or~(\ref{EqGOne}) and~(\ref{EqGTwo})]
on the priors $g_S$ for the nonzero coefficients ensures
that the posterior does not shrink to the center
of the prior too much. In the next theorem we investigate the necessity of
conditions of this type. The theorem shows
that product priors with marginal densities proportional to $y\mapsto
e^{-|y|^\alpha}$ for some $\alpha>1$ lead to a slow contraction rate
for large true vectors $\q_0$. We formulate this in an asymptotic
setting with a sequence of true vectors, written as $\q_0^n$,
tending to infinity. We denote by $p_n$ the number
of nonzero coordinates of $\q_0^n$.

The theorem applies in particular to the normal distribution.
For this prior a problem (only) arises if the parameter
vector $\q_0^n$ tends to infinity
faster than the optimal rate
\[
\bigl\|\te_0^n\bigr\|^2\gg p_n
\log(n/p_n).
\]
The posterior then puts no mass on balls of radius a multiple of
$\|\q_0^n\|$ around the true parameter. For ``small'' $\q_0^n$ no
problem occurs, because shrinkage to the origin is desirable in that
case. However, if the true parameter satisfies $\|\te_0^n\|^2\lesssim
p_n\log(n/p_n)$, then the estimator that is zero, irrespective of the
observations, possesses mean square error of the order the minimax
risk for the problem. Thus it is rather poor consolation that the
Bayes procedure based on Gaussian priors performs well in this case,
as it is no better than the ``zero estimator.'' Gaussian priors really
are problematic.

Product priors with marginal density proportional to
$y\mapsto e^{-|y|^a}$ give behavior as the Gaussian prior for every
$\alpha\ge2$.
For $\alpha\in(1,2)$ the result is slightly more complicated and involves
the quantities
%
%e2.11 #&#
%
\begin{equation}
\label{rna} \rho_{0,\al}^n = \biggl(\frac{\|\te_0^n\|_{\alpha}^{\al}}{\|
\te_0^n\|_2^2}\wedge
1 \biggr) \bigl\|\te_0^n\bigr\|_{\al}
p_n^{1/2-1/\al},
\end{equation}
where \mbox{$\|\cdot\|_{\al}$} denotes the usual $L_{\al}$-norm on $\RR^n$
(i.e.,
$\|\q\|_\alpha^\alpha=\sum_i|\q_i|^\alpha$). The numbers $\rho_{0,\alpha
}^n$ increase to infinity
as $\q_0^n$ tends to infinity at a sufficiently fast rate. For instance
$\rho_{0,\alpha}^n$ is of the order $c_n^{\alpha
-1}p_n^{1/2-1/\alpha}$ if $\alpha<2$
and $\q_0^n=c_n\bar\q_0$
for scalars $c_n$ and fixed vectors $\bar\q_0$. The following
theorem shows that if $\rho_{0,\alpha}^n$ increases to infinity
faster than
the optimal rate $ (p_n\log(n/p_n) )^{1/2}$, then the posterior
does not charge balls of radius a small multiple of $\rho_{0,\alpha}^n$.
%
%th2.8 #&#
%
\begin{theorem}[(Heavy tails)]
\label{TheoremLB}
Assume that the densities $g_S$ are products of $S$ univariate
densities proportional to $y\mapsto e^{-|y|^\alpha}$ and the
prior $\pi_n$ on dimension satisfies~(\ref{EqLBOnPin}) for some $c>0$:
%$$\pi_n(p_n) \ge\exp(-cp_n\log(n/p_n)).$$
%
\begin{longlist}[(ii)]
\item[(i)] If $\alpha\ge2$ and $ \|\te_0^n\|^2/ (p_n\log
(n/p_n) )
\to\infty$, then for sufficiently small $\h>0$, as $n\to\infty$,
\[
P_{n,\q_0^n}\Pi_n \bigl(\q\dvtx\bigl\|\te-\te_0^n
\bigr\|\leq\h\bigl\|\te_0^n\bigr\| \given X^n\bigr) \ra0.
\]
\item[(ii)] If $1<\alpha<2$ and $(\rho_{0,\al}^n)^2/(p_n\log
(n/p_n))\ra\infty$, then
for sufficiently small $\h>0$, as $n\to\infty$,
\[
P_{n,\q_0^n}\Pi_n \bigl(\q\dvtx\bigl\|\te-\te_0^n
\bigr\|\leq\h\rho_{0,\al
}^n\given X^n\bigr)\ra0.
\]
\end{longlist}
\end{theorem}

Theorem~\ref{TheoremLB} shows problematic behavior of
the posterior distribution for signals with large
energies $\|\te_0^n\|$. Instead of using fixed priors on the
coordinates, we could make them depend on the sample size, for
instance, Gaussian priors with variance $v_n\to\infty$, or
uniform priors on intervals $[-K_n,K_n]$ with $K_n\to\infty$.
Such priors will push the ``problematic boundary'' toward
infinity, but the same reasoning as for the theorem will show
that shrinkage remains for (very) large $\q_0^n$.

The above results show that $g_S$ needs to have heavy
tails. Another important condition, this time concerning the prior
$\pi_n$ on the dimension $k$, concerns the amount of mass $\pi_n(p_n)$
at the true dimension. If this quantity is too small, then the Bayes
procedure might not be optimal.
%
%th2.9 #&#
%
\begin{theorem}
\label{tlbdim}
Suppose also that the prior $\pi_n$
on dimension in \textup{(P1)} is decreasing and that there exist integers
$d_{1,n}<d_{2,n}$ such that, for some $C>0$ and a sequence $\epl$ such
that $n\epl^2\to\infty$,
\[
\frac{\pi_n(d_{2,n})}{\pi_n(d_{1,n})} \pmatrix{n \cr d_{1,n} } \le
e^{-Cn\epl^2}.
\]
Denoting $d_{3,n}=(3d_{2,n}-d_{1,n})/2$, there exists $\te_0$ in $\ell
_0[d_{3,n}]$ such that,
for sufficiently small $\h>0$, as $n\to\infty$,
\[
P_{n,\q_0^n}\Pi_n \bigl(\q\dvtx\bigl\|\te-\te_0^n
\bigr\|\leq\h\rn\epl\given X^n\bigr) \ra0.
\]
\end{theorem}
%
%({\bf Prior on dimension in $e^{-k\log{k}}$}).
%If $\pi_n(k)=r\exp(-k\log{k})$, with $r$ the appropriate
%normalizing constant, let us apply the preceding result
%with the choices $d_{1,n}=p_n/4$, $d_{2,n}=3p_n/4$,
%for some sequence $p_n\to\infty$. It holds
%$$ \frac{\pi_n(3p_n/4)}{\pi_n(p_n/4)}{ n \choose p_n/4 }
%e^{-\frac{p_n}{4}\log{p_n}-\frac{p_n}{4}
%If we assume that $p_n^2/n\to\infty$, the last display is smaller
%than $\exp(-Cp_n\log p_n)$ for $C=1/4$. Thus Theorem~\ref{tlbdim}
%implies
% that there is a vector $\te_0$ in $\ell[p_n]$ such that
%$$P_{n,\q_0^n}\Pi_n\bigl(\q: \|\te-\te_0^n\|^2\leq\h p_n\log p_n
%for a small enough constant $\h$.
%This implies that the corresponding estimator does not reach the
%minimax rate over the class
%$\ell_0[p_n]$ as soon as $p_n\log p_n$ tends to infinity faster than
%$p_n\log(n/p_n)$
%(think for instance of $p_n$ of the form $n/\exp(\sqrt{\log{n}})$).
%
%ex2.7 #&#
%
\begin{example}[{[Prior on dimension in $\exp(-k(\log{k})^a)$, with
$a\ge1$]}]
If $\pi_n(k)=r\exp(-k\log^a{k})$, with $r$ the appropriate
normalizing constant, let us apply the preceding result
with the choices\vadjust{\goodbreak} $d_{1,n}=p_n/4$, $d_{2,n}=3p_n/4$, for some sequence
$p_n\to\infty$. It holds
\begin{eqnarray*}
\frac{\pi_n(3p_n/4)}{\pi_n(p_n/4)}\pmatrix{n \cr p_n/4 } & \le&
e^{-({3p_n}/{4})\log^a ({3p_n}/{4} )
+({p_n}/{4})\log^a ({p_n}/{4} )+({p_n}/{4})\log
(ne)}
\\
& \le&e^{-({p_n}/{4})\log^a ({3p_n}/{4} )
-({p_n}/{4})\log^a ({3p_n}/{4} )^{2^{1/a}}+
({p_n}/{4})\log^a(ne)}.
\end{eqnarray*}
As long as we impose $(3p_n/4)^{2^{1/a}}\geq ne$ and
$\log(3p_n/4)\geq2^{-1/a}\log{p_n}$ (which holds for large enough $n$),
the last display is at most $\exp(-\frac{p_n}{8}\log^a{p_n})$.
Theorem~\ref{tlbdim} implies that there is a vector $\te_0$ in $\ell
_0[p_n]$ with
\[
P_{n,\q_0^n}\Pi_n \bigl( \q\dvtx\bigl\|\te-\te_0^n
\bigr\|^2 \leq\h p_n\log^a p_n \given
X^n \bigr) \ra0
\]
for a small enough constant $\h$.
This implies that the corresponding estimator does not reach the
optimal rate over the class
$\ell_0[p_n]$ as soon as $p_n\log^a p_n$ tends to infinity faster
than $p_n\log(n/p_n)$ [take, e.g., $p_n=n/\exp(\sqrt{\log{n}})$].
\end{example}

%s2.7 #&#
\subsection{Discussion}

We have identified general conditions on the prior that ensure optimal
convergence rates for estimating a sparse mean vector in Gaussian
noise. In particular, natural fully Bayes estimates
(e.g., Beta-binomial prior on dimension) are shown to be adaptive with
respect to the unknown smoothing parameter $p_n/n$.

% General idea of: downweight big models and use heavy tails.
Especially in high-dimensional contexts the full
posterior measure and special aspects of it can start to have
divergent behaviors. We have seen that for nonconvex distances the
posterior mean is not a satisfactory projection. It can also happen
that the mode and the full posterior behave differently.

%Lower bounds What if underpenalize, find examples
In some situations one might want to estimate prior hyperparameters,
and in this case, it is desirable to assess the convergence
properties of the resulting plug-ins. To our knowledge, there are only
a few works in this direction; see~\cite{JohnSilv04,jiangzhang}.
Potential alternative proofs could consist in
obtaining first (suitably uniform) results for the (full) posterior
measure and combine them with a statement saying that ``the plug-in
estimate is not too bad.'' Also, here, one could evaluate the sparsity
coefficient $\eta_n=p_n/n$ via the posterior number $\hat{k}_n$ of
selected models and plug this estimate into the full posterior for the
binominal prior on dimension. Since $\hat{\eta}_n=\hat{k}_n/n$ does
not exceed $Cp_n$ with high probability, we have some control of the
plug-in into the full posterior. The question of then deriving results
for estimates of it (e.g., the mean), remains open.

%s3 #&#
\section{Implementation}
\label{SectionSimulations}
In this section we provide an algorithm to compute several functionals
of the posterior measure associated with the prior defined by \mbox{(P1)--(P3)},
including the posterior mean, marginal posterior quantiles
and the posterior of the number of selected models.
The algorithm is exact in that it does not rely on an
approximation of the posterior distribution, but computes the
exact expressions. We illustrate the
posterior quantities through simulations.

We assume that the densities $g_S$ on $\RR^S$ are products of $S$
copies of a univariate density $g$. Because the prior on the number of
nonzero coordinates induces dependence, this generally does not entail
a factorization of the posterior distribution as a product measure.
(An exception is the binomial distribution for $\pi_n$.)

For all computations, we need the denominator of the posterior measure
in~(\ref{EqPosterior}) (the ``partition function'').
For $\phi$ the standard normal
density, and $\psi=\phi\ast g$ its convolution with the density $g$,
this can
be written
\[
Q_n:=\sum_{p=0}^{n}
\frac{\pi_n(p)}{{n\choose p}} \sum_{|S|=p} \prod
_{i\in S}\psi(X_i)\prod
_{i\notin S}\phi(X_i).
\]
Naive computation directly from this expression would
require a number of operations that grows exponentially with $n$.
However, the sum over all models $S$ of size~$p$
(the inner sum in the display) is
equal to the coefficient of $Z^p$ in the polynomial
\[
Z\mapsto\prod_{i=1}^n \bigl(
\phi(X_i)+ \psi(X_i) Z\bigr).
\]
This polynomial can be computed by a quadratic number of operations
by computing the products term by term, and in $n\log^2{n}$ operations
by a more clever algorithm.

%s3.1 #&#
\subsection{Posterior mean}
The posterior mean $\hat{\te}^{PM}=\int\q\,d\Pi_n(\q\given X)$
is a random vector in $\RR^n$. Letting
$\zeta(x) = \int t \phi(x-t) g(t) \,dt$,
we can write its first coordinate in the form
\[
\hat\q^{PM}_1=\frac{1}{Q_n} \sum
_{p=1}^{n} \frac{\pi_n(p)}{{n\choose p}} \zeta(X_1)
\sum_{|S|=p, 1\in S} \prod_{i\in S, i\neq1}
\psi(X_i) \prod_{i\notin S}
\phi(X_i).
\]
The inner sum (over $S$) is the coefficient of $Z^p$ in the polynomial
$Z\mapsto\zeta(X_1)Z \prod_{i=2}^n (\phi(X_i)+ \psi(X_i) Z)$.
Hence it can be computed as before.

%s3.2 #&#
\subsection{Coordinatewise quantiles}
The distribution function of the marginal posterior distribution
of the first coordinate can be written, for any real $u$,
\[
\Pi\bigl( (-\infty,u]\times\RR^{n-1} \given X \bigr) =
(1-q_{n,1}) 1_{u\ge0} + q_{n,1} \frac{\psi(X_1,u)}{\psi(X_1)},
\]
where $1-q_{1,n}$ is the posterior probability that the first coordinate
is zero, and $\psi(x, u) = \int_{-\infty}^{u} \phi(x-t) g(t) \,dt$.
The former probability can be written
\[
1-q_{n,1}=\Pr(\q_1=0\given X) =\frac1{Q_n}
\sum_{p=0}^{n} \frac{\pi_n(p)}{{n\choose p}} \sum
_{|S|=p, 1\notin S} \prod_{i\in S}
\psi(X_i) \prod_{i\notin S}
\phi(X_i).
\]
Hence it can be computed as before, now involving the polynomial
$Z\mapsto\psi(X_1)Z \prod_{i=2}^n (\phi(X_i)+ \psi(X_i) Z)$.

Given the marginal posterior distribution, we can compute
marginal quantiles. For instance, the first component of the
coordinatewise median $\hat{\te}^{\mathrm{med}}$ is given by,
with $H_{n,1}^{-1}$ the inverse of $H_{n,1}(u)=\psi(X_1,u)/\psi(X_1)$,
\[
\hat{\te}^{\mathrm{med}}_1 = \biggl[H_{1,n}^{-1}
\biggl(1-\frac
{1}{2q_{1,n}} \biggr) \vee0 \biggr] + \biggl[H_{n,1}^{-1}
\biggl(\frac{1}{2q_{n,1}} \biggr)\wedge0 \biggr].
\]
The last display should be understood with the convention
$H_{n,1}^{-1}(u)=-\infty$ if $u\le0$ and $H_{n,1}^{-1}(u)=\infty$ if
$u\ge1$.
%It can be explicitely written as
%$$ 1-q_{1,n} = \frac{ \phi(X_1)}{Q_n}{ {\sum_{p=0}^{n}

%s3.3 #&#
\subsection{Number of nonzero coordinates}
The posterior distribution of the number
$|S_\q|$ of nonzero coordinates of $\q\in\RR^n$ is the random distribution
on the set $\{0,1,\ldots, n\}$ given by
\[
\Pi_n \bigl(\q\dvtx|S_\q|=p\given X \bigr) =\frac1{Q_n}
\frac{\pi_n(p)}{{n\choose p}} \sum_{|S|=p} \prod
_{i\in S}\psi(X_i) \prod
_{i\notin S} \phi(X_i).
\]
The same computational scheme applies. In fact the sum
will already be computed in the derivation of $Q_n$.

%s3.4 #&#
\subsection{Simulations}
In a small simulation study we considered
the prior defined by (P1)--(P3) with $g$ a Laplace density
$x\to(a/2)e^{-a|x|}$, with scale parameter \mbox{$a>0$} and
two priors on dimension, suggested by our theoretical results,
given by
%
%e3.1 #&#
%e3.2 #&#
%
\begin{eqnarray}
\label{EqPrior1} \pi_n(p)&\propto& e^{-\kappa p\log(3n/p) },
\\
\label{EqPrior2} \pi_n(p)&\propto& \pmatrix{2n-p \cr n}^{\kappa}.
\end{eqnarray}
Here $\kappa$ is a real parameter, which for both priors quantifies
how fast they decrease to zero with $p$.
In the results shown we used $a=1$ and $\kappa\in\{0.1,1\}$.

We simulated signals $\q=(\q_1,\ldots,\q_n)$ of length $n=500$,
for various settings of the sparsity $p_n=\#(\q_i\not=0)$
and for signals $\q$ with the nonzero coordinates set equal to
a fixed number $A$. We show the results for $p_n\in\{25,50,100\}$
and ``signal strength'' $A\in\{3, 4, 5\}$.

%t1 #&#
%
\begin{table}
\caption{Average square errors of eight estimators computed on 100
data vectors $X$ of length $n=500$ simulated from model
(\protect\ref{EqModel}) with $\q=(0,0,\ldots,0,A,\ldots,A)$, where
$p_n$ coordinates indices are equal to $A$.
In every column the smallest value is printed in bold face.
The estimators are: \textup{PM1}, \textup{PM2}: posterior means for
two priors $\pi_n$ in
(\protect\ref{EqPrior1}) and (\protect\ref{EqPrior2})
and Laplace prior on nonzero coordinates;
\textup{PMed1}, \textup{PMed2} coordinatewise medians for the same priors;
\textup{EBM}, \textup{EBMed}: empirical Bayes mean
and median for Laplace prior; \textup{HT}, \textup{HTO}:
hard-thresholding and
hard-thresholding-oracle}
\label{table1}
\begin{tabular*}{\tablewidth}{@{\extracolsep{\fill}}l c r c c c c c c c@{}}
\hline
$\bolds{p_n}$ & \multicolumn{3}{c}{\textbf{25}} &
\multicolumn{3}{c}{\textbf{50}} &  \multicolumn{3}{c@{}}{\textbf{100}}
\\[-4pt]
& \multicolumn{3}{c}{\hrulefill} &
\multicolumn{3}{c}{\hrulefill} &  \multicolumn{3}{c@{}}{\hrulefill} \\
$\bolds{A}$ & \textbf{3} & \multicolumn{1}{c}{\textbf{4}} & \textbf{5} &  \textbf{3} & \textbf{4} & \textbf{5}&
\textbf{3} & \textbf{4} & \textbf{5} \\
\hline
\\     [-4pt]
PM1&111&96& 94& 176 &165 &154& 267&302&307\\
PM2& \textbf{106}& 92& 82 & 169 &165& 152& 269&280&274\\
EBM &\textbf{103}& 96& 93& 166& 177& 174& 271&312&319 \\
PMed1 &129& \textbf{83}& 73 & 205& 149 & \textbf{130} & \textbf{255}&279 &283
\\
PMed2 &125& 86& \textbf{68} & 187&148& \textbf{129}& 273&\textbf{254}
&\textbf{245} \\
EBMed & 110 & \textbf{81}& 72 & \textbf{162}& 148&142&\textbf{255}&294&300 \\
HT&175& 142& \textbf{70}& 339& 284& 135 & 676&564&252 \\
HTO &136& 92& 84 & 206& 159& 139&306&261 &245 \\
\hline
\end{tabular*}
\end{table}

%t2 #&#
%
\begin{table}[b]
\caption{Average absolute deviation errors of eight estimators
computed on 100
data vectors $X$ of~length $n=500$ simulated from model
(\protect\ref{EqModel}) with $\q=(0,0,\ldots,0,A,\ldots,A)$, where
$p_n$ coordinates indices are equal to $A$. In every column the
smallest value
is printed in bold face. The priors and estimators are
as in Table \protect\ref{table1}}
\label{table2}
\begin{tabular*}{\tablewidth}{@{\extracolsep{\fill}}l c r r r r r c c c@{}}
\hline
$\bolds{p_n}$ & \multicolumn{3}{c}{\textbf{25}} &
\multicolumn{3}{c}{\textbf{50}} &  \multicolumn{3}{c@{}}{\textbf{100}}
\\[-4pt]
& \multicolumn{3}{c}{\hrulefill} &
\multicolumn{3}{c}{\hrulefill} &  \multicolumn{3}{c@{}}{\hrulefill} \\
$\bolds{A}$ & \textbf{3} & \multicolumn{1}{c}{\textbf{4}}
& \multicolumn{1}{c}{\textbf{5}} &  \multicolumn{1}{c}{\textbf{3}}
& \multicolumn{1}{c}{\textbf{4}} & \multicolumn{1}{c}{\textbf{5}} &
\multicolumn{1}{c}{\textbf{3}} & \multicolumn{1}{c}{\textbf{4}}
& \multicolumn{1}{c@{}}{\textbf{5}} \\
\hline
\\    [-4pt]
PM1 &80 & 101& 110 & 127 & 145 &147& 240&268&270\\
PM2& 79 & 85& 87 & 135 &145& 144& 219&232&232\\
EBM & 95& 110& 117& 191& 200& 176& 260&285&281 \\
PMed1 &\textbf{51}& 43& 45& \textbf{86} & \textbf{80} & 78 & 178&225&230 \\
PMed2 & \textbf{50} & \textbf{40}& 37 & \textbf{86}& \textbf{79}& 76&
\textbf{156}&\textbf{162}&163 \\
EBMed & \textbf{50} & 48& 45 & 108& 121& 97& 212 & 258 & 257 \\
HT& 63& 44& \textbf{27}& 122& 86& \textbf{53} & 244&173&\textbf{102} \\
HTO& 53& 41& 40 & 91& 79& 74& 157&148&144 \\
\hline
\end{tabular*}
\end{table}

Tables~\ref{table1} and~\ref{table2} report estimates of
the mean square errors $\E\|\hat\q-\q\|^2$ and mean
absolute deviation errors $\E\|\hat\q-\q\|_1$ of eight estimators
$\hat\q$. These estimates are the average
(square) error of 100 estimates $\hat\q_1,\ldots,\hat\q_{100}$
computed from $100$ data vectors simulated
independently from model~(\ref{EqModel}). The eight estimators include the
posterior means \textit{PM1}, \textit{PM2} and coordinatewise medians
\textit{PMed1}, \textit{PMed2} associated with the two priors $\pi_n$ with
$\kappa=0.1$,
the empirical Bayes mean \textit{EBM} and median \textit{EBMed}
considered in~\cite{JohnSilv04} with a standard Laplace prior,
and the hard-thresholding \textit{HT} and hard-thresholding-oracle
\textit{HTO} estimators, given by
\[
\hat\q^{\mathrm{HT}}_i=X_i1_{|X_i|>\sqrt{2\log{n}}},\qquad \hat
\q^{\mathrm{HTO}}_i=X_i1_{|X_i|>\sqrt{2\log{n/p_n}}}.
\]
The last estimator uses the ``oracle'' value of the sparsity
parameter $p_n$, whereas the other seven estimators do not use
this value.\vadjust{\goodbreak}

The tables show that the mean and median of the full Bayesian
posterior distribution are competitive with the empirical Bayes
estimates. The behavior of the full Bayes and empirical Bayes
estimates seems similar, up to a few aspects. In terms of
squared risk, empirical Bayes estimates appear to be slightly better for
small $p_n$ and small $A$, while the full Bayes estimates appear to
be slightly better for larger signals and larger $p_n$. For
$L^1$-risk, the full Bayes estimates appear to outperform the
EB-estimates in
most of the cases. (Additional simulation results, not shown,
suggest that the situation becomes less unfavorable
for empirical Bayes as the scale parameter $a$ of the Laplace prior
is taken smaller than $1$.) In agreement
with~\cite{JohnSilv04}, in most cases the mean estimates perform not
quite as well as the median ones, already in terms of squared-risk.

The parameter $a$ of the Laplace prior plays the same role for
the full Bayes as for the empirical Bayes estimates.
Although we do not investigate this aspect here, it could
be estimated from the data, as is proposed in the
EbayesThresh package, or be treated as a hyperparameter in a full
Bayes approach. [A single scale parameter for high-dimensional
densities $g_S$ appears to create dependence between the coordinates
that is stronger
than what is allowed by our conditions~(\ref{EqGOne}) and~(\ref{EqGTwo}), and
hence would need further analysis.]
Similar remarks pertain to the parameter $\kappa$. The
choice $\kappa=0.1$ seemed to be fairly good uniformly over all
considered simulations, also for smaller $n$'s.
%Taking $\kappa_1$ of about $1$ or larger makes in a way the prior go
%to zero too fast which results in a decrease of performance of the
%estimates.

\begin{figure}[t!]

\includegraphics{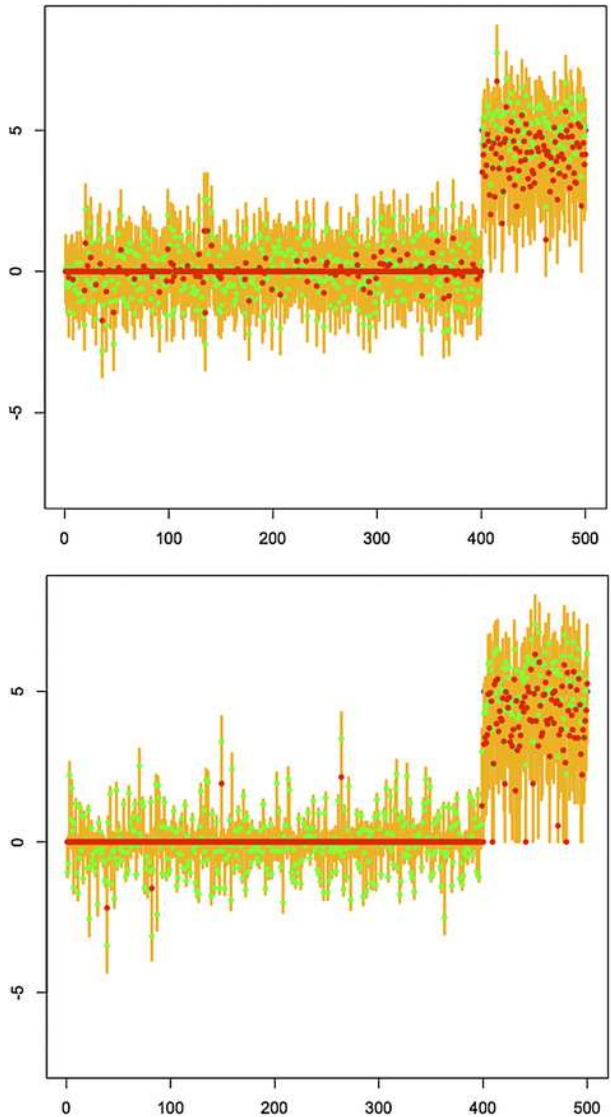}

\caption{Marginal posterior medians (red dots) and marginal credible intervals
(orange) for the parameters $\q_1,\ldots,\q_n$ for a single
data vector $X_1,\ldots,X_n$ simulated according to the model
(\protect\ref{EqModel}) with $\q=(0,0,\ldots,0,5,\ldots,5)$, where
$n=500$ and the last $p_n=100$ coordinates are nonzero. The
data points are indicated by green dots.
The prior $g$ is the standard Laplace density, and $\pi_n$ is as in
(\protect\ref{EqPrior2}) with ``inverse temperature'' $\k_1=0.1$
(TOP graph)
and $\k_1=1$ (BOTTOM graph).}
\label{Figure1}
\end{figure}

For further illustration Figure~\ref{Figure1} shows marginal 95\%
credible intervals (orange bars) for the parameters
$\q_1,\ldots,\q_n$, and marginal posterior medians (red dots) for a
single simulation of the data vector, with single strength $A=5$,
$p_n=100$ and $n=500$.
The observations
$X_1,\ldots,X_n$ are indicated by green dots. The credible intervals
are defined as intervals between the 2.5\% and 97.5\% percent
quantiles of the marginal posterior distributions of the
parameters. The intervals corresponding to zero and nonzero
coefficients $\q_i$ are clearly separated, although some of the
credible intervals of nonzero $\q_i$ contain the value zero. Also
visible is that the posterior medians and the credible intervals
surrounding them are shrunk toward zero
relative to the observed value $X_i$, for the zero coordinates $\q_i$,
which is desirable, but also for the nonzero $\q_1$.
Figure~\ref{Figure1} (bottom) shows that for $\k=1$ the shrinkage
effects are stronger, and the credible intervals shorter.

%f1 #&#
%

Since our main goal here is illustration, we
only implemented a simple version of the algorithm. This computes
the polynomials with direct loops and can be improved. This
implementation is limited to $n$ of the order
500, not by computing time, but by the appearance of large
numbers in the polynomial coefficients that overflow standard
memory capacity $(10^{-300},10^{300})$. Handling
larger $n$ should certainly be possible by improved programming,
for instance, by computing on a logarithmic scale.
Algorithmic complexity appears not to be a major issue.

%
%intervals
%(orange) for the parameters $\q_1,\ldots,\q_n$ for a single
%data vector $X_1,\ldots,X_n$ simulated according to the model
%$n=500$ and the last $p_n=100$ coordinates are nonzero. The
%data points are indicated by green dots.
%The prior $\pi_n$ is as in \eqref{EqPrior2}
%with ``inverse temperature'' $\k_1=1$,
%and $g$ is the standard Laplace density.}

%s4 #&#
\section{\texorpdfstring{Proof of Theorem \protect\ref{TheoremDimension}}{Proof of Theorem 2.1}}
We first prove the theorem for priors on dimension $\pi_n(p)$ with
strict exponential decrease and densities $g_S$ that are not
necessarily of
product form, but that satisfy
(\ref{EqBoundConditionalDensity}),
for $Dm_1<1$, and $D$ the constant in~(\ref{priordim}). Thus the
proof also covers half of Theorem~\ref{TheoremConcentrationDependent}.
In view of Example~\ref{exproduct}, densities of the product form
% of a univariate density
%with mean zero and finite second moment
satisfy~(\ref{EqBoundConditionalDensity})
with $m_1=1$, and hence automatically have $Dm_1<1$.

Since the true parameter $\q_0$ is assumed to have $p_n$ nonzero
coordinates, it is sufficient to prove that the intersection of the
support $S_\q$ with the complement $S_0^c$ of the support
$S_0\triangleq S_{\q_0}$ of $\q_0$ has dimension of the order $p_n$
under the posterior distribution. The following proposition gives an explicit
bound on this dimension; it is followed by a lemma that shows
that this bound tends to zero under the conditions of the theorems.
The idea of the proof of the proposition is to
condition on the vector of the coordinates $\pi_{S_0}\q$ of $\q$ that
belong to $S_0$.

The unconditional density of
$(S_\q,\q)$ for $\q$ drawn from the prior $\Pi_n$ is given by, with
$\delta_0$
denoting a ``Dirac density at $0$,''
\[
(S,\q)\mapsto\frac{\pi_n(|S|)}{{n\choose|S|}}g_S(\q_S)
\delta_0(\q_{S^c}).
\]
The conditional density of $(S_\q\cap S_0^c,\q_{S_0^c})$
given $\q_{S_0}$ is proportional
to this expression viewed as function of $(S\cap S_0^c,\q_{S\cap S_0^c})$.
This\vspace*{-1pt} shows the conditional distribution has the same structure
as the prior $\Pi_n$, but with sample space $\RR^{S_0^c}$ rather than
$\RR^n$,
with the density of the nonzero coordinates of $\q_{S_0^c}$ given
by $g_{S\cap S_0^c|S\cap S_0}(\cdot\given\q_{S\cap S_0})$,
proportional to
$g_{S\cap S_0^c,S\cap S_0}(\cdot,\q_{S\cap S_0})$,
and the prior on dimension given by
%
%e4.1 #&#
%
\begin{equation}
\label{EqDefPink} \pi_{n,k}(p)\propto\pi_n(p+k)
\frac{{n-p_n\choose p}}{{n\choose
p+k}},\qquad k=|S_\q\cap S_0|.
\end{equation}
The extra factor (quotient) on the right arises because
$\pi_{n,k}(p)$ and $\pi_n(p+k)$ are the probabilities of the given
dimensions, and hence the sums of the probabilities of all subsets of
that dimension. Recall also that we assume that $\pi_n(p)$ is positive
for any $p$,
which makes the maximum appearing in the following proposition always finite.
%
%pr4.1 #&#
%
\begin{proposition}
\label{PropositionDimension}
If the densities $g_S$ satisfy~(\ref{EqBoundConditionalDensity}),
%Theorem~\ref{TheoremDimension}or Theorem
then, for any $ A\ge1$,
%, with $2\nu_k \triangleq m_2 \sum_{p=0}^{n-p_n} p\pi_{n,k}(p)$,
%
\[
\sup_{\q_0\in\ell_0[p_n]}P_{n,\q_0}\Pi_n\bigl(\q\dvtx
\bigl|S_\q\cap S_{\q
_0}^c\bigr|\ge A \given X \bigr) \leq
\sum_{p=A}^{n-p_n} m_1^{p_n+p}
\max_{0 \le k \le p_n} \biggl[ \frac{ \pi_{n,k}(p) }{ \pi_{n,k}(0)
} \biggr].
\]
\end{proposition}
\begin{pf}
For $B= \{\q\dvtx|S_\q\cap S_{0}^c|\ge A \}$ and
$\Pi_n^{\q_{S_0}}(\cdot\given X)$ the marginal distribution of $\q
_{S_0}$ if $\q$ is
distributed according to the posterior distribution,
\begin{eqnarray*}
\Pi_n(B\given X) & =&\int\Pi_n(B\given X,
\q_{S_0}=\bar\te_1) \,d\Pi_n^{\q_{S_0}}(\bar
\te_1\given X)
\\
& \leq& \sup_{\bar\te_1\in\RR^{S_0}} \Pi_n(B\given X,\q_{S_0}=\bar
\te_1).
\end{eqnarray*}
In the Bayesian\vspace*{1pt} setting the vectors $X_{S_0}$ and $X_{S_0^c}$ are
conditionally independent given $\q$ with marginal conditional
distributions depending on $\q_{S_0}$ and $\q_{S_0^c}$ only, respectively.
This implies that the distribution of $\q_{S_0^c}$ given $(X,\q_{S_0})$
depends\vspace*{1pt} on $(X_{S_0^c},\q_{S_0})$ only. The joint distribution
of $(X_{S_0^c},\q_{S_0^c},\q_{S_0})$ can be\vspace*{1pt} generated by first generating
$\q_{S_0}$ from its marginal distribution derived from $\Pi_n$, next
generating
$\q_{S_0^c}$ from its conditional given $\q_{S_0}$ derived from $\Pi_n$,
and finally generating $X_{S_0^c}$ from the $N_{n-p_n}(\q
_{S_0^c},I)$-distribution.
It\vspace*{1pt} follows that the conditional distribution
of $\q_{S_0^c}$ given $(X,\q_{S_0})$ can also be described
as the ``ordinary'' posterior distribution of $\q_{S_0}^c$ given
the observation $X_{S_0^c}$ relative to the prior on $\q_{S_0^c}$
given by the conditional distribution of $\q_{S_0^c}$ given $\q_{S_0}$ derived
from $\Pi_n$. If $\Pi_n(\cdot\given\bar\te_1)$ denotes
the prior induced on $\RR^{S_0^c}$ when conditioning $\Pi_n$ to the event
that $\q_{S_0}=\bar\q_1$, and $\bar n_2=n-p_n$, then
%
%e4.2 #&#
%
\begin{equation}
\label{EqHULPDIM} \Pi_n(B\given X,\q_{S_0}=\bar
\te_1) = \frac{\int_{B} p_{\bar n_2,\bar\te_2}(X_{S_0^c}) \,d\Pi_n(\bar
\te_2\given\bar\q_1)} {
\int p_{\bar n_2,\bar\te_2}(X_{S_0^c}) \,d\Pi_n(\bar\te_2\given
\bar\q_1)}.
\end{equation}
%
%The denominator of the right side can be bounded
%below with the help of Lemma~\ref{lemlowbound} (below),
%applied with $\bar n_2$ instead of $n$,
%with both $\Pi$ and $\tilde\Pi$ taken equal to $\Pi_n(\cdot\given\bar
%and with $\q_0=0$. This yields the bound for the right side of
%the preceding display given by
%$$ e^{\sigma_2^2/2 -\mu_2^TX_{S_0^c}} \int_{B}
% d\Pi_n(\bar\te_2\given\bar\q_1), $$
%where $\m_2=\int\bar\q_2 d\Pi_n(\bar\q_2\given\bar\q_1)$
%and $\s_2^2=\int\|\bar\q_2\|^2 d\Pi_n(\bar\q_2\given\bar\q_1)$.
%In fact $\m_2=0$, by the assumption that the conditional distributions
%corresponding to the $g_S$ have zero mean. If the nonzero
%coordinates of $\bar\q_1$ have index in $S_1\subset S_0$
%and values $\q_1\in\RR^{S_1}$, then
%for $|S_1|=k$ and $\Pi_n(S_2\given S_1)$ the probability that
%a vector $\q$ drawn from the prior has $S_\q\cap S_0^c=S_2$ given that
%$S_\q\cap S_0=S_1$,
% & \leq& \sum_{S_2|S_0^c} \Pi_n(S_2\given S_1)m_2 |S_2|
The denominator of the right-hand side can be bounded below by
restricting the integrating set to the singleton $\{\bar\te_2=0\}$,
leading to
\[
\int p_{\bar n_2,\bar\te_2}(X_{S_0^c}) \,d\Pi_n(\bar
\te_2\given\bar\q_1) \ge\Pi_n(\bar
\q_2=0 \given\bar\q_1) p_{\bar n_2,0_{S_0^c}}(X_{S_0^c}).
\]

Let $S_2$ denote the indices of the nonzero coordinates of $\bar\q_2\in
\RR^{S_0^c}$,
$\q_2$ the vector of their values and $n_2=|S_2|$, and similarly for
$S_1,\te_1$. Then
%$$\frac{p_{\bar n_2,\bar\q_2}}{p_{\bar n_2,0_{S_0^c}}}(X_{S_0^c})
%=\frac{p_{n_2,\q_2}}{p_{n_2,0_{S_2}}}(X_{S_2}).$$
%Therefore
%
\begin{eqnarray*}
\Pi_n(B\given X, \q_{S_0}=\bar\q_1) &\le&
\Pi_n(\bar\q_2=0 \given\bar\q_1)^{-1}
\int_B \frac{p_{\bar
n_2,\bar\q_2}}{p_{\bar n_2,0_{S_0^c}}}(X_{S_0^c}) \,d
\Pi_n(\bar\q_2\given\bar\q_1)
\\
&\le& \sum_{S_2\subset S_0^c,|S_2|\ge A} \frac{\Pi_n(S_2\given\bar\q
_1)}{ \Pi_n(S_2=\varnothing\given\bar
\q_1)} \int
\frac{p_{\bar n_2,\bar\q_2}}{p_{\bar
n_2,0_{S_0^c}}}(X_{S_0^c}) \,d\Pi_n(\bar\q_2
\given\bar\q_1,S_2). %& \le\sum_{p=A}^{n-p_n}
% d\q_{S_2}.
\end{eqnarray*}
With the notation $S_1,\te_1,\te_2$ introduced above, one obtains
\[
\int\frac{p_{\bar n_2,\bar\q_2}}{p_{\bar
n_2,0_{S_0^c}}}(X_{S_0^c}) \,d\Pi_n(\bar
\q_2\given\bar\q_1,S_2) = \int
\frac{p_{n_2,\te_2}}{p_{n_2,0_{S_2}}}(X_{S_2}) \frac{
g_{S_1,S_2}(\te_1,\te_2)}{
\int g_{S_1,S_2}(\te_1,\te_2) \,d\te_2}\,d\te_2.
\]
On the other hand, an application of Bayes's formula leads to
\[
\frac{\Pi_n(S_2\given\bar\q_1)}{ \Pi_n(S_2=\varnothing\given\bar
\q_1)} = \frac{\Pi_n(S_1,S_2)}{ \Pi_n(S_1,S_2=\varnothing)} \int\frac
{g_{S_1,S_2}(\te_1,\te_2)}{g_{S_1}(\te_1)}\,d
\te_2,
\]
and the last ratio of prior probabilities of subsets is equal to
\[
\frac{\Pi_n(S_1,S_2)}{ \Pi_n(S_1,S_2=\varnothing)} = \frac{\pi_n(p+k)}{{n
\choose p+k}} \frac{ {n \choose k} }{\pi_n(k)}=\frac{\pi_{n,k}(p)} {
\pi_{n,k}(0)}
\frac{1}{{n-p_n \choose p}}.
\]
Combining the previous identities and condition
(\ref{EqBoundConditionalDensity}), one obtains that $\Pi_n(B\given
X,
\q_{S_0}=\bar\q_1)$ is bounded above, uniformly in
$\bar\q_1\leftrightarrow(S_1,\q_1)$, by
\[
\sum_{p=A}^{n-p_n} \sum
_{|S_2|=p} \max_{0\le k\le p_n} \biggl[ \frac{ \pi_{n,k}(p) }{ \pi
_{n,k}(0) } \biggr]
\frac{m_1^{p_n+p}}{{n-p_n\choose p}} \int\frac{p_{n_2,\te
_2}}{p_{n_2,0_{S_2}}}(X_{S_2})\ga_{S_2}(
\te_2) \,d\te_2.
\]
The proposition follows,
since $P_{n,\te_0}p_{n_2,\te_2}/p_{n_2,0_{S_2}}(X_{S_2})=1$.
\end{pf}
%
%le4.1 #&#
%
\begin{lemma} \label{lemexpde2}
%If $\pi_n(p)\le D \pi_n(p-1)$ for some constant $D<1$ and all $p\ge
%Cp_n$,
%If $\pi_n$ satisfies \eqref{priordim}, then $\nu_k\le m_2 D_1 p_n$
%for a constant $D_1$ that depends on $C, D$ in \eqref{priordim} only.
If $\pi_n$ satisfies~(\ref{priordim}) with $C=0$ and a constant $D$
such that $m_1D<1$, then $\sum_{p= P_n}^{n-p_n}
m_1^{p_n+p}\max_k [\pi_{n,k}(p)/\pi_{n,k}(0)] \ra0$ for
$P_n$ bigger than a sufficiently large multiple of $p_n$ and $P_n\ra
\infty$.
\end{lemma}
\begin{pf}
From the expression of $\pi_{n,k}$ in~(\ref{EqDefPink}), simple
algebra leads to
\[
\frac{\pi_{n,k}(p)}{\pi_{n,k}(0)} = \pmatrix{p+k \cr k} \frac{\pi
_n(p+k)}{\pi_n(k)} \frac{(n-p_n)\times\cdots\times
(n-p_n-p+1) } {
(n-k)\times\cdots\times(n-k-p+1) }.
\]
Using the assumed strict exponential decrease, the second ratio in the
last display is bounded above by $e^{p\log{D}}$. For any integer $k$
between $0$ and $p_n$, the last factor (ratio) in the last display is
bounded above by $1$ and ${p+k \choose k}$ is bounded above by ${p+p_n
\choose p_n}\le e^{p_n\log\{e (p+p_n)/p_n\}}$. Since $\log(1+x)\le
x/M$, for $M>0$ as soon as $x$ is larger than a sufficiently large
multiple of $M$, the result follows.
%The numbers $b_p:={n-p_n\choose p}/{n\choose p+k}$ satisfy
%$$\frac{b_p}{b_{p-1}}=\frac{n-p_n-p+1}{n-p-k+1}\frac{p+k}p
%=\Bigl(1-\frac{p_n-k}{n-p-k+1}\Bigr)\Bigl(1+\frac kp\Bigr).$$
%Hence $b_p\le b_{p-1}(1+C_1^{-1})$ for $p\ge C_1p_n$ and any $k\le
%p_n$, whence
%the numbers $a_p=\pi_n(p+k) b_p$ satisfy
%$$a_p\le a_{p-1} D(1+C_1^{-1})\le a_{C_1p_n} F^{p-C_1p_n},$$
%for $p\ge C_1p_n$ and $F=D(1+C_1^{-1})$, provided $C_1$ is larger than
%the constant $C$ in the assumptions.
%Because $D<1$ and $m_1D<1$, there exists suitable $C_1$ such that
%$F<1$ and $m_1 F<1$ (for the next line we use only $F<1$). Then
%$$\frac{\sum_{p=C_1p_n}^{n-p_n}p a_p}{\sum_{p=0}^{n-p_n}a_p}
%F^{p-C_1p_n}}{a_{C_1p_n}}
%The missing initial part of the normalized sum in the left side also
%contributes at most $C_1p_n$. This concludes the proof
%for the bound on
%%the expectations
%$\nu_k$.
%
%For the final assertion we must take the dependence of $a_p=a_{p,k}$
%on $k$ into account. The preceding argument shows
%that $a_{p,k}\le a_{Cp_n,k} F^{p-Cp_n}$ for constants $C$ and $F$ that
%do
%not depend on $k$. Therefore $a_{p,k}/\sum_p a_{p,k}\le F^{p-Cp_n}$
%for every
%$p\ge Cp_n$ and every $k$. It follows that the sum in the lemma
%is bounded above by $\sum_{p\ge P_n} m_1^p F^{p-Cp_n} e^{m_2 D_1 p_n}$,
%for $P_n\ge Cp_n$, where $m_1F<1$.
\end{pf}

Combining Proposition~\ref{PropositionDimension} and Lemma \ref
{lemexpde2} concludes the proof of the first half of Theorem \ref
{TheoremConcentrationDependent} and of Theorem~\ref{TheoremDimension} for priors on
dimension with strict exponential decrease.

For $g_S$ of the product form and $\pi_n$ with just exponential
decrease [$C>0$ in~(\ref{priordim})] such as the oracle binomial
prior, we use a slight variant of the above argument. Starting from
(\ref{EqHULPDIM}),
the denominator can be bounded
below with the help of Lemma~\ref{lemlowbound} (below),
applied with $\bar n_2$ instead of $n$,
with $\q_0=0$ and both $\Pi=\tilde\Pi=\Pi_n(\cdot\given\bar\te_1)$.
This implies
that $\Pi_n(B\given X,\te_{S_0}=\bar\te_1)$
is bounded above by
\[
e^{\sigma_2^2/2 -\mu_2^TX_{S_0^c}} \int_{B} \frac{p_{\bar n_2,\bar\te
_2}}{p_{\bar n_2,0_{S_0^c}}}(X_{S_0^c})
\,d\Pi_n(\bar\te_2\given\bar\q_1),
\]
where $\m_2=\int\bar\q_2 \,d\Pi_n(\bar\q_2\given\bar\q_1)$
and $\s_2^2=\int\|\bar\q_2\|^2 \,d\Pi_n(\bar\q_2\given\bar\q_1)$.
In fact $\m_2=0$, by the assumption that the common density $g$ has
zero mean. If $m_2$ denotes the second moment of $g$, we have
\[
\s_2^2 = \sum_{S_2|S_0^c}
\Pi_n(S_2\given S_1) m_2
|S_2| \le m_2\sum_{p=0}^{n-p_n}
p \pi_{n,k}(p) \triangleq2 \nu_k.
\]
This implies that $\Pi_n(B\given X,\te_{S_0}=\bar\te_1)$ is
uniformly bounded in $\bar{\te}_1$ by
\[
\sum_{p=A}^{n-p_n} \sum
_{|S_2|=p} \max_{0\le k\le p_n} \bigl(\pi_{n,k}(p)e^{\nu_k}
\bigr) \frac{1}{{n-p_n\choose p}} \int\frac{p_{n_2,\te
_2}}{p_{n_2,0_{S_2}}}(X_{S_2})g_{S_2}(
\te_2) \,d\te_2.
\]
To conclude one takes the $P_{n,\te_0}$-expectation and uses Lemma
\ref{lemexpde} below.
%
%le4.2 #&#
%
\begin{lemma} \label{lemexpde}
%If $\pi_n(p)\le D \pi_n(p-1)$ for some constant $D<1$ and all $p\ge
%Cp_n$,
If $\pi_n$ satisfies~(\ref{priordim}), then $\nu_k\le m_2 D_1 p_n$
with $D_1$ that depends on $C, D$ in~(\ref{priordim}) only.
Furthermore, $\sum_{p= P_n}^{n-p_n}
\max_k (\pi_{n,k}(p)e^{\nu_k})\ra0$ for
$P_n$ bigger than a sufficiently large multiple of $p_n$ and $P_n\ra
\infty$.
\end{lemma}

%The numbers $b_p:={n-p_n\choose p}/{n\choose p+k}$ satisfy
%$$\frac{b_p}{b_{p-1}}=\frac{n-p_n-p+1}{n-p-k+1}\frac{p+k}p
%=\Bigl(1-\frac{p_n-k}{n-p-k+1}\Bigr)\Bigl(1+\frac kp\Bigr).$$
%Hence $b_p\le b_{p-1}(1+C_1^{-1})$ for $p\ge C_1p_n$ and any $k\le
%p_n$, whence
%the numbers $a_p=\pi_n(p+k) b_p$ satisfy
%$$a_p\le a_{p-1} D(1+C_1^{-1})\le a_{C_1p_n} F^{p-C_1p_n},$$
%for $p\ge C_1p_n$ and $F=D(1+C_1^{-1})$, provided $C_1$ is larger than
%the constant $C$ in the assumptions.
%Because $D<1$ and $m_1D<1$, there exists suitable $C_1$ such that
%$F<1$ and $m_1 F<1$ (for the next line we use only $F<1$). Then
%$$\frac{\sum_{p=C_1p_n}^{n-p_n}p a_p}{\sum_{p=0}^{n-p_n}a_p}
%F^{p-C_1p_n}}{a_{C_1p_n}}
%The missing initial part of the normalized sum in the left side also
%contributes at most $C_1p_n$. This concludes the proof
%for the bound on
%%the expectations
%$\nu_k$.
%
%For the final assertion we must take the dependence of $a_p=a_{p,k}$
%on $k$ into account. The preceding argument shows
%that $a_{p,k}\le a_{Cp_n,k} F^{p-Cp_n}$ for constants $C$ and $F$ that
%do
%not depend on $k$. Therefore $a_{p,k}/\sum_p a_{p,k}\le F^{p-Cp_n}$
%for every
%$p\ge Cp_n$ and every $k$. It follows that the sum in the lemma
%is bounded above by $\sum_{p\ge P_n} m_1^p F^{p-Cp_n} e^{m_2 D_1 p_n}$,
%for $P_n\ge Cp_n$, where $m_1F<1$.
%

%%%%%%%%%%%%%%%%%%%%%%%%%%
%%%% Proof in Supplementary file %%%%
%%%%%%%%%%%%%%%%%%%%%%%%%%

%s5 #&#
\section{\texorpdfstring{Proof of Theorems \protect\ref{TheoremConcentrationIndependent} and \protect\ref{TheoremConcentrationDependent}}
{Proof of Theorems 2.2 and 2.4}}
\label{secproof2}
\label{secmeanmed}
In view of Theorem~\ref{TheoremDimension} the posterior mass of
models of dimension bigger than $Ap_n$, for a large constant $A$,
tends to zero. Thus it suffices to show concentration
around $\q_0$ in models with $|S_\q|\le Ap_n$. This is achieved
using testing arguments. Proposition~\ref{Theorem2} gives an
explicit bound on concentration with respect to the Euclidean
metric. General $d_q$-metrics are next treated by interpolation
of metrics.

Let $\Phi$ be the standard normal distribution function
and $\bar\Phi=1-\Phi$.
%
%le5.1 #&#
%
\begin{lemma}
\label{LemmaTestPhi}
For any $\alpha,\b>0$ and any $\q_0,\q_1\in\RR^n$ there exists a test
$\phi$ based on $X\sim N(\q,I)$, such that for every
$\q\in\RR^n$ with $\|\q-\q_1\|\le\|\q_0-\q_1\|/2\triangleq\rho$,
\[
\alpha P_{n,\q_0}\phi+\b P_{n,\q}(1-\phi) \le\alpha\bar{\Phi}
\biggl(\frac\rho2+\frac1\rho\log\frac\alpha\b\biggr) +\b\Phi\biggl
(-\frac\rho2+
\frac1\rho\log\frac\alpha\b\biggr).
\]
This quantity can be further bounded by
$2\sqrt{\al\be} e^{-\|\te_0-\te_1\|^2/32}$.
%The function on the right is increasing in $\b$, decreasing
%in $\rho$, and bounded above by $\b$.
\end{lemma}

%%%%%%%%%%%%%%%%%%%%%%%%%%
%%%% Proof in Supplementary file %%%%
%%%%%%%%%%%%%%%%%%%%%%%%%%

We note that the bound of Lemma~\ref{LemmaTestPhi}, even though valid
for every $\alpha,\b>0$, is of interest only if $\alpha$ and $\b$
are not too
different: if $\log\alpha/\b\le-\|\q_0-\q_1\|^2/32$ or $\log
\alpha/\b\ge
\|\q_0-\q_1\|^2/32$, then the trivial tests $\phi=1$ and $\phi=0$ give
the better bounds $\alpha$ and $\b$, respectively.
%
%le5.2 #&#
%
\begin{lemma} \label{lemlowbound}
For any prior probability distribution $\Pi$ on $\RR^n$,
any positive measure $\tilde\Pi$ with $\tilde\Pi\le\Pi$,
and any $\q_0\in\RR^m$,
\[
\int\frac{p_{n,\q}}{p_{n,\q_0}}(X) \,d\Pi(\q) \ge\|\tilde\Pi\|
e^{-\tilde
\s^2/2+\tilde\m^T(X-\q_0) },
\]
where $\tilde\m=\int(\q-\q_0) \,d\tilde\Pi(\q)/\|\tilde\Pi\|$ and
$\tilde\s^2=\int\|\q-\q_0\|^2 \,d\tilde\Pi(\q)/\|\tilde\Pi\|$.
Consequently, for any $r>0$,
\[
P_{n,\q_0} \biggl(\int\frac{p_{n,\te}}{p_{n,\te_0}} \,d\Pi(\te) \geq e^{-r^2}
\Pi\bigl(\q\dvtx\|\q-\q_0\|<r \bigr) \biggr)\geq1 - e^{-r^2/8}.
\]
\end{lemma}
%
%le5.3 #&#
%
\begin{lemma}
\label{LemmaVolume}
The volume $v_p$ of the $p$-dimensional Euclidean unit ball satisfies,
for every $p\ge1$, setting $d_1=1/\sqrt{\pi}$ and $d_2=e^{1/6}d_1$,
\[
d_1 (2e\pi)^{p/2}p^{-p/2-1/2} \leq v_p
\leq d_2 (2e\pi)^{p/2} p^{-p/2-1/2}.
\]
\end{lemma}
%
%%%%%%%%%%%%%%%%%%%%%%%%%%
%%%% Proof in Supplementary file %%%%
%%%%%%%%%%%%%%%%%%%%%%%%%%
%
%le5.4 #&#
%
\begin{lemma} \label{lembeta}
Let $S\subset\{1,\ldots,n\}$, $p=|S|$, $j\ge1$ and $r_n^2 \geq p_n
\vee\log\pi_n(p_n)^{-1}$. Let $\te_{S,j}\in\RR^n$ with support
$S$ and $2j r_n< \| \te_{S,j}-\te_0 \|<2(j+1)r_n$. For some universal
constant $c_3>0$, we have that
%Assume \eqref{EqGOne}-\eqref{EqGTwo} on the prior densities and that
%$r_n^2\ge p_n\log(n/p_n)\vee\log\bigl(1/\pi_n(p_n)\bigr)$. Then
%
\begin{eqnarray*}
&& \log{\frac{\Pi(\te\in\RR^n\dvtx S_\te=S, \|\pi_S\te-\te_{S,j}\|
<jr_n)}{e^{-r_n^2}
\Pi(\te\in\RR^n, \|\te-\te_0\|<r_n)} }
\\
&&\qquad\le c_3(p + p_n) + p\log{j}+ 9(j+1)^2
r_n^2/64 + 7 r_n^2 /2.
\end{eqnarray*}
\end{lemma}
\begin{pf}
Denoting $\b_{S,j}$ the quantity in the logarithm in the last display,
\begin{eqnarray*}
\b_{S,j} %&=\frac{\Pi(\te\in\RR^S\dvtx\|\te-\te_{S,j}
&\le&\frac{\Pi(S)G_S (\q\in\RR^S\dvtx\|\q-\pi_S\q_{S,j}\|
<jr_n )} {
e^{-r_n^2}\Pi(S_0)G_{S_0} (\q\in\RR^{S_0}\dvtx\|\q-\pi_{S_0}\q_0\|<r_n )}
\\
&\le&\frac{\Pi(S)v_{S}(jr_n)^{|S|}\max(g_{S}(\q)\dvtx
\|\q-\pi_{S}\q_{S,j}\|<jr_n ) } {
e^{-r_n^2}\Pi(S_0)v_{S_0}r_n^{|S_0|}
\min(g_{S_0}(\q)\dvtx\|\q-\pi_{S_0}\q_0\|<r_n )}.
\end{eqnarray*}
Let us decompose, for any $\q'\in\RR^{S}$ and $\q\in\RR^{S_0}$,
\[
\frac{g_{S}(\q')}{g_{S_0}(\q)} = \frac{g_{S}(\q')}{g_{S\cap S_0}(\pi
_{S\cap S_0}\q')} \frac{g_{S\cap S_0}(\pi_{S\cap S_0}\q')}{g_{S\cap
S_0}(\pi_{S\cap
S_0}\q)} \frac{g_{S\cap S_0}(\pi_{S\cap S_0}\q)}{g_{S_0}(\q)}.
\]
Combining this identity with~(\ref{EqGOne}) and~(\ref{EqGTwo}), we
obtain, with
$c_2=1/64$,
\begin{eqnarray*}
\biggl| \log\frac{g_{S}(\q')}{g_{S_0}(\q)} \biggr| &\le& c_1|S|+c_1|S \cap
S_0|+c_1|S_0|
\\
&&{} +c_2 \bigl\| \pi_{S-S_0} \q' \bigr\|^2
+c_2 \bigl\| \pi_{S\cap S_0} \bigl(\q'-\q\bigr)
\bigr\|^2 +c_2 \| \pi_{S_0-S} \q\|^2.
\end{eqnarray*}
Denoting by $\bar{\te},\bar{\te}'$ the vectors of $\RR^n$ with
respective supports $S_0, S$ and such that $\pi_{S_0}\bar{\te}=\te
$, $\pi_S\bar{\te'}=\te'$, note that the last line of the previous
display is bounded above by $c_2\|\bar{\te}'-\bar{\te}\|^2$.
For $\|\q'-\pi_{S}\q_{S,j}\|<jr_n$ and $\|\q-\pi_{S_0}\q_0\|<r_n$,
we have %(note that $\pi_{S_0}\q_{S,j}=(\pi_{S\cap S_0}
\[
\bigl\|\bar{\te}'-\bar{\te}\bigr\| \le\bigl\| \bar{\te}'-
\te_{S,j}\bigr\| +\|\te_{S,j}-\te_0\| +\|
\te_0-\bar{\te}\| \le3(j+1)r_n. %\lefteqn{ \| \pi_{S-S_0} \q' \|^2 +\|
%+ \| \pi_{S_0-S} \q\|^2}\\
%&&\le(\| \pi_{S-S_0} \q' \|^2 +\| \pi_{S\cap S_0} (\q'-\q) \|^2 )
%+ (\| \pi_{S\cap S_0} (\q'-\q) \|^2 + \| \pi_{S_0-S} \q\|^2 )\\
%&&\le jr_n+2(j+1)r_n+r_n=3(j+1)r_n.
\]
Due to Lemma~\ref{LemmaVolume}, the quotient $v_p
r_n^p/(v_{p_n}r_n^{p_n})$ is bounded by
\[
\frac{v_p r_n^{p/2}}{v_{p_n} r_n^{p_n}} \lesssim(2e\pi)^{p} \biggl(\frac
{\sqrt{p_n}}{r_n}
\biggr)^{p_n} \biggl(\frac{r_n}{\sqrt{p}} \biggr)^{p}.
\]
Since $r_n^2\ge p_n$ by assumption, we have $(\sqrt{p_n}/r_n)^{p_n}\le
1$, and because the function $p\mapsto p\log(r_n^2/p)$ takes a maximum
at $p=r_n^2/e$, we obtain, for some universal constants $C,C'$,
\[
\b_{S,j}\le j^p e^{Cp+C'p_n+9c_2(j+1)^2r_n^2 + (1+1/2e)r_n^2}\Pi(S)/
\Pi(S_0). %&\le e^{c_3(p+p_n) + p\log{j}+ 9c_2 (j+1)^2 r_n^2 + (3+1/2e)
%r_n^2 }.
\]
To conclude, one notes that $\Pi(S)\le1$ and
that ${n\choose p_n}\le(ne/p_n)^{p_n}\le e^{r_n^2+p_n}$ by the
assumption on $r_n$,
so that $\Pi(S_0)\ge e^{-2r_n^2-p_n}$.
%and $9(j+1)^2/64+3+\log j\le j^2/4$ for $j\ge10$.
\end{pf}
%
%pr5.1 #&#
%
\begin{proposition}
\label{Theorem2}
If the densities $g_S$ satisfy~(\ref{EqGOne}) and~(\ref{EqGTwo})
and have finite second moments, then there exist
universal constants $d_1,d_2, d_3$ such that
for \mbox{$M\ge10$} and $1\le A\le n/(2p_n)$ and
$r_n^2$ satisfying~(\ref{EqConditionRate}) and $p_n/n\to0$,
as $n\to\pli$,
\begin{eqnarray*}
&&
\sup_{\te_0\in\ell_0[p_n]} P_{n,\te_0} \Pi_n\bigl(\te\dvtx\|\te-
\te_0\|>M r_n, |S_\q|\le Ap_n
\given X\bigr)
\\
&&\qquad\le e^{-r_n^2/8} + d_1 \pmatrix{n \cr A p_n }
e^{d_2Ap_n -d_3(Mr_n)^2}.
\end{eqnarray*}
\end{proposition}
\begin{pf}
Let $\cS_1$ be the collection of subsets $S\subset\{1,2,\ldots,n\}$
such that $|S|\le Ap_n$.
For each such $S$ and $j=1,2,\ldots$
let $\{\q_{S, j,i}\dvtx i\in I_{S,j}\}$ be a maximal $jr_n$-separated
set inside the set $\{\q\in\RR^n\dvtx S_\q=S, 2jr_n\le\|\q-\q_0\|
\le2(j+1)r_n\}$.
Because the latter set is within
a ball of radius $2(j+1)r_n$ of the projection $\Pi_S\q_0$ onto
the subspace of vectors with support\vspace*{1pt} inside $S$,
a volume argument shows that
the cardinality of $I_{S,j}$ is at most $9^{|S|}$.

We can partition the set of vectors with exactly support $S$ by assigning
each such vector to a closest point $\q_{S,j,i}$ for some
$j=1,2,\ldots,$
and $i\in I_{S,j}$. The resulting partitioning sets $B_{S,j,i}$
will fit into balls of radius $jr_n$. For each
$\q_{S,j,i}$ fix a test $\phi_{S,j,i}$ as in Lemma~\ref{LemmaTestPhi}
with $\alpha=1$ and the triple $(\q_0,\q_1)$, $\rho$ and $\b$
taken equal to
the triple $(\q_0,\q_{S,j,i})$, $jr$ and $\b_{S,j,i}$, where the
last numbers
will be determined later.
In view of the second assertion of Lemma~\ref{lemlowbound} applied
with $r$ equal to $r_n$, there exist events $\A_n$ such that
$P_{n,\te_0}(\A_n^c)\le e^{-r_n^2/8}$, on which
\[
\int\frac{p_{n,\te}}{p_{n,\te_0}} \,d\Pi_n(\te) \geq e^{-r_n^2}
\Pi_n \bigl(\q\dvtx\|\q-\q_0\|< r_n \bigr).
\]
We have that
\begin{eqnarray*}
&&
P_{n,\q_0} \Pi_n \bigl(\q\dvtx\|\q-\q_0\|>
2Mr_n, S_\q\in\cS_1\given X
\bigr)1_{\A_n}
\\
&&\qquad\le\sum_{S\in\cS_1}\sum_{j\ge M}
\sum_{i \in I_{S,j}} P_{n,\q_0}\Pi_n (\q\in
B_{S,j,i}\given X )1_{\A_n}
\\
&&\qquad\le\sum_{S\in\cS_1}\sum_{j\ge M}
\sum_{i \in I_{S,j}} \biggl(P_{n,\q_0}
\phi_{S,j,i}\\
&&\qquad\quad\hspace*{69.2pt}{}+P_{n,\q_0} \biggl[(1-\phi_{S,j,i})
\frac{\int_{B_{S,j,i}}{p_{n,\q}}/{p_{n,\q
_0}} \,d\Pi(\q)} {
e^{-r_n^2}\Pi(\q\dvtx\|\q-\q_0\|<r_n )} \biggr] \biggr)
\\
&&\qquad\le\sum_{S\in\cS_1}\sum_{j\ge M}
\sum_{i \in I_{S,j}} \Bigl(P_{n,\q_0}
\phi_{S,j,i}+\b_{S,j,i}\sup_{\q\in
B_{S,j,i}}P_{n,\q}(1-
\phi_{S,j,i}) \Bigr),
\end{eqnarray*}
where we have denoted
\[
\b_{S,j,i}=\frac{\Pi(B_{S,j,i})}{e^{-r_n^2}\Pi(\q\dvtx\|\q-\q_0\|<r_n )}.
\]
In view of Lemma~\ref{LemmaTestPhi} the term within the triple sum is
bounded using by $2\sqrt{\be_{S,j,i}}e^{-j^2r_n^2/8}$. Since
$|S|=p\le
Ap_n$ and $p_n/n\to0$, we can take $n$ large enough in order to
have both $c_3(p+p_n) \le r_n^2/10$ and $p\log j\le j^2r_n^2/100$ for
any\vspace*{1pt} $j\ge1$. Since $M\ge10$, we have $j\ge10$, so we also have
$r_n^2\le j^2r_n^2/100$.

Combination
with Lemma~\ref{lembeta} now yields the bound, for $j\ge10$,
\[
\log\sqrt{\b_{S,j,i}} \le2.3 j^2r_n^2/100+9(j+1)^2r_n^2/128.
\]
One easily checks that this is bounded by $(1-d_2)j^2r_n^2/8$, for
$d_2=1/9$ when $j\ge10$. Thus the probability at stake is bounded from
above by
\[
%& \sum_{S\in\cS_1}\sum_{j\ge M}\sum_{i \in I_{S,j}} e^{-c_3j^2r_n^2}\\
\sum_{p=0}^{Ap_n} \pmatrix{n \cr p}
\sum_{j\ge M} 2 C^{p} e^{-d_2j^2r_n^2} \leq
d_1^{Ap_n}e^{-d_2 M^2 r_n^2} \sum
_{p=0}^{Ap_n} \pmatrix{n \cr p}
\]
for $d_1$ large enough. By assumption $Ap_n\le n/2$, so
each binomial term is bounded by the last one. Using simple algebra
this yields the second term in the bound of the theorem. The first
term comes from $P_{n,\te_0}1_{\cA_n^c}\le e^{-r_n^2/8}$.~%
\end{pf}

In view of~(\ref{EqConditionRate}) we have
${n\choose Ap_n}\le(ne/Ap_n)^{Ap_n} \le e^{d_4r_n^2}$. Therefore, the
right-hand
side of Proposition~\ref{Theorem2} tends to zero.
Combining this with Theorem~\ref{TheoremDimension}
yields proofs of Theorems~\ref{TheoremConcentrationIndependent}
and~\ref{TheoremConcentrationDependent}
for $d_q$ the square Euclidean norm $d_2$.

The theorems for $q\in(0,2)$ are a corollary of the case $q=2$,
by interpolation between the distances.
Due to H\"older's inequality, for any $\q,\q_0$ with $|S_\q\cup
S_0|\le Ap_n$,
\[
d_q(\te,\te_0) \leq\|\te-\te_0
\|^q (Ap_n)^{1-q/2}.
\]
This implies, for any $M>0$, if $\q_0\in\ell_0[p_n]$,
\begin{eqnarray*}
&&P_{n,\te_0} \Pi_n \bigl(d_q(\te,
\te_0)>M r_n^qp_n^{1-q/2}
\given X \bigr) \\
&&\qquad\le P_{n,\te_0} \Pi_n\bigl(\q\dvtx
|S_\q|>(A-1)p_n\given X\bigr)
\\
&&\qquad\quad{}+ P_{\te_0}^n \Pi\bigl( \|\te-\te_0\| >
M^{1/q}A^{1/2-1/q} r_n \given X \bigr).
\end{eqnarray*}
The first term on the right-hand side tends to zero for sufficiently large
$A$. Next the second tends to zero for sufficiently large $M$.

%%%%%%%%%%%%%%%%%%%%%%%%%%%%%%%%%%%%%%%%%%
%%%% Section "Proofs for Complexity priors" in Supplementary file %%%%
%%%%%%%%%%%%%%%%%%%%%%%%%%%%%%%%%%%%%%%%%%

%s6 #&#
\section{\texorpdfstring{Proof of Theorem \protect\ref{TheoremWeakClass}}{Proof of Theorem 2.6}}
The theorem is proved by bounding the (posterior) risk under
a vector $\q_0\in m_s[p_n]$ by the risk under its projection
into $\ell_0[p]$ obtained by setting the smallest
$n-p$ coordinates of $\q_0$ equal to zero. The value $p$
that minimizes the expression that defines the rate $r_n^2$ is
the optimal dimension of a projection, and the complicated
expression itself is a trade-off of an approximation error
and a rate.\vadjust{\goodbreak}

The comparison between $\q_0$ and its projection $\q_1$ is made
in the following lemma.
%
%le6.1 #&#
%
\begin{lemma}\label{postcomp}
For any measurable function $G$ and any
$\te_0,\te_1$ in $\mathbb{R}^n$,
\[
P_{n,\q_0}G \leq\sqrt{P_{n,\te_1}G^2}
e^{\|\te_1-\te_0\|^2/2}.
\]
\end{lemma}
\begin{pf}
In view of the Cauchy--Schwarz inequality,
\[
P_{n,\q_0}G \leq\sqrt{ P_{n,\te_1} G^2} \sqrt{
P_{n,\q_1} \biggl(\frac{dP_{n,\q_0}}{dP_{n,\te_1}} \biggr)^2 }.
\]
The second integral on the right-hand side is equal to $\exp(\|
\q_0-\q_1\|^2 )$.
\end{pf}

Let $p_n^*$ be an index for which the minimum that defines
the rate $r_n^2$ is attained. For given $\te_0$ belonging to $m_s[p_n]$,
let $\te_1$ denote the vector deduced from $\te_0$ by keeping unchanged
its $p_n^*$ largest components and putting the other ones to $0$.
By definition $\te_1$ belongs to $\ell_0[p_n^*]$ and
%
%e6.1 #&#
%
\begin{eqnarray}\label{EqDistanceTheta01}
\|\te_1 - \te_0\|^2 & = & \sum
_{ i > p_n^* } |\te_{0,[i]}|^2 \leq
\biggl(\frac{p_n}{n} \biggr)^2 \sum_{ i > p_n^* }
\biggl(\frac
{n}{i} \biggr)^{2/s}
\nonumber\\[-8pt]\\[-8pt]
& \leq& \biggl(\frac{p_n}{n} \biggr)^2 \biggl(
\frac{s}{2-s} \biggr) n^{2/s} \bigl(p_n^*
\bigr)^{1-2/s} \le r_n^2,\nonumber
\end{eqnarray}
where the first inequality is obtained
using the definition of the $m_s[p_n]$-class, and the second
follows by comparison of the series with an integral.

Therefore, the triangle inequality implies
\[
\Pi_n \bigl(\q\dvtx\|\q-\q_0\|> 80r_n+20r \given X
\bigr) \leq\Pi_n \bigl(\q\dvtx\|\q-\q_1\|> 79r_n +20r
\given X \bigr).
\]
By Lemma~\ref{postcomp} the expectation of the right-hand side
under $P_{n,\q_0}$ is bounded by
\[
\bigl(P_{n,\q_1}\Pi_n \bigl(\q\dvtx|\q-\q_1
\|>79r_n+20r \given X \bigr) \bigr)^{1/2} e^{\|\q_0-\q_1\|^2/2}.
\]
Finally apply
Theorem~\ref{TheoremMainExponentialDecrease}, with $r$ of the theorem taken
equal to $3.4r_n+2r$.

%s7 #&#
\section{\texorpdfstring{Proof of Theorems \protect\ref{TheoremLB} and \protect\ref{tlbdim}}
{Proof of Theorems 2.8 and 2.9}}
\label{secproofslb}
The proof of Theorem~\ref{TheoremLB}
follows the approach to get lower bound type results
introduced in~\cite{ic08}, which uses the principle
that sets with very little prior mass
receive no posterior mass, see also Figure~\ref{fig2}.
%
%le7.1 #&#
%
\begin{lemma}
\label{lemlb}
We have $P_{n,\q_0}\Pi_n (\q\dvtx\|\q-\q_0\|<s_n\given X
)\ra0$,
for any $s_n$ for which there exist $r_n$ such that
\[
\frac{\Pi_n (\q\dvtx\|\q-\q_0\|< s_n )} {
\Pi_n (\q\dvtx\|\q-\q_0\|< r_n )}=o\bigl(e^{-r_n^2}\bigr).
\]
\end{lemma}
%
%le7.2 #&#
%
\begin{lemma}\label{lemvol}
There exist a constant $C>0$ such that if $S\subset\{1,\ldots,n\}$
and $r_n$ is a sequence of real numbers such that $r_n^2\ge|S_{\q
_0}|$, it holds
\[
\frac{v_{|S \cap S_{\q_0}|}}{v_{|S_{\q_0}|}} \frac{1}{r_n^{|S_{\q
_0}\setminus S|}} \leq e^{C |S_{\q_0}|}.
\]
\end{lemma}
\begin{pf*}{Proof of Theorem~\ref{TheoremLB}}
We first consider the (more complicated) case that $1<\alpha<2$.
For this range of $\alpha$ an application of
H\"older's inequality gives that $\|\q\|_\alpha\le\|\q\|
p^{1/\alpha-1/2}$,
%
%f2 #&#
%
\begin{figure}

\includegraphics{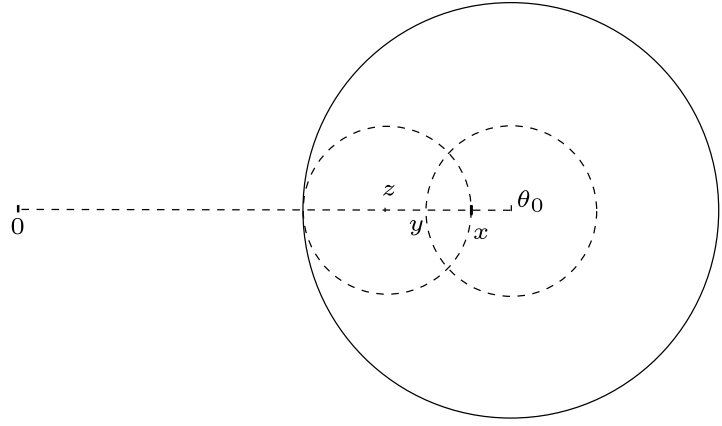}

\caption{Idea behind the proof of Theorem \protect\ref{TheoremLB}.}\label{fig2}
\end{figure}
if $p$ is the number of nonzero coordinates of a vector $\q$.
Let us introduce
\[
r_n= \biggl(\frac{\|\te_0\|_\al^\al}{\|\te_0\|^2}\wedge1 \biggr) \frac
{\|\te_0\|}{8},\qquad
s_n = \frac{\rho_{0,\alpha}^n}{64}=\frac{r_n}{8} \biggl(
\frac{\|\te_0\|_\al}{\|\te_0\|}p_n^{{1}/{2}-
{1}/{\al}} \biggr).
\]
Then $r_n \leq\|\te_0\|/8$ and $s_n\leq r_n/8$. Also, %Let
%$\Pi_n(S)=\pi_n(p)/{n\choose p}$ denote the prior probability
%of selecting a parameter with support a particular set
%$S\subset\{1,\ldots, n\}$ of
%cardinality $p$, and let $S_0=S_{\q_0}$ be the support of $\q_0$. Then
%
\begin{eqnarray*}
&&
\frac{\Pi_n (\q\dvtx\|\q-\q_0\|< s_n )} {
\Pi_n (\q\dvtx\|\q-\q_0\|< r_n )}
\\
&&\qquad= \sum_S \Pi_n(S)\frac{G_S (\q\in\RR^S\dvtx
\|\q-\pi_S\q_0\|^2+\|\pi_{S_0\setminus S}\q_0\|^2< s_n^2 )} {
\Pi_n (\q\dvtx\|\q-\q_0\|< r_n )}
\\
&&\qquad\leq\sum_S\frac{\Pi_n(S)}{\Pi_n(S_0)} \frac{G_{S\cap S_0} (\q
\in\RR
^{S\cap S_0}\dvtx
\|\te-\pi_{S\cap S_0}\te_0\|\le s_n )} {
G_{S_0} (\q\in\RR^{S_0}\dvtx\|\te-\pi_{S_0}\te_0\|\le r_n )}
1_{\|\pi_{S_0\setminus S}\q_0\|<s_n}.
\end{eqnarray*}
%
%Here we have used that the marginal onto $\RR^{S\cap S_0}$
%of the product prior $G_S$ coincides with $G_{S\cap S_0}$.
Define
\[
\q_B= \biggl(1-\frac{r_n-s_n}{\|\te_0\|} \biggr)\pi_{S_0}
\te_0^n.
\]
Then the ball in $\RR^{S_0}$ of radius $s_n$ around $\q_B$ is contained
in the ball of radius $r_n$ around $\pi_{S_0}\q_0$. It follows that
the second-to-last display is bounded above by
%
%e7.1 #&#
%
\begin{equation}
\label{EqQuotientProofLB} \sum_S
\frac{\Pi_n(S)}{\Pi_n(S_0)} \frac{s_n^{|S\cap S_0|} v_{S\cap
S_0}}{s_n^{p_n} v_{p_n}} \frac{{\sup_{\q\in A}g_{S\cap S_0}(\q)}} {
{\inf_{\q\in B}g_{S_0}(\q)}} 1_{ \llVert\pi_{S_0\setminus S}\te
_0\rrVert\le s_n }
\end{equation}
with $A=\{\q\in\RR^{S\cap S_0}\dvtx\| \q-\pi_{S\cap S_0}\te_0^n \|
<s_n\}$ and
$B=\{\q\in\RR^{S_0}\dvtx\|\q-\q_B\|<s_n\}$.
We finish the proof by bounding the densities $g_{S\cap S_0}$
and $g_{S_0}$ above and below on the given sets.

If $\q\in B$, then by the triangle inequality
followed by H\"older's inequality,
\begin{eqnarray*}
\|\q\|_\al& \leq& \|\te_B\|_\al+ \|\q-
\te_B\|_\al
\\
& \leq& \biggl(1-\frac{r_n-s_n}{\|\te_0\|} \biggr) \|\te_0
\|_\al+ p_n^{{1}/{\al}-{1}/{2}} s_n \le\biggl(1-
\frac{3r_n}{4\|\te_0\|} \biggr) \|\te_0\|_\al,
\end{eqnarray*}
because $s_n\le r_n/8$ and
$p_n^{{1}/{\al}-{1}/{2}} s_n\le(r_n/8)\|\q_0\|_\alpha/\|
\q_0\|$.
Similarly, if $\q\in A$ and $\|\pi_{S_0\setminus S}\q_0\|<s_n$, then
$\|\pi_{S_0\setminus S}\q_0\|_\alpha<p_n^{1/\alpha-1/2} s_n$ and
\begin{eqnarray*}
\|\q\|_{\al} &\geq& \|\q_0\|_\alpha-\|
\q_0-\pi_{S\cap S_0}\q_0\|_\al- \|
\pi_{S\cap S_0}\q_0-\q\|_\alpha
\\
& \geq& \|\te_0\|_{\al} - 2p_n^{{1}/{\al}-{1}/{2}}
s_n \geq\|\te_0\|_{\al} \biggl(1-
\frac{r_n}{4\|\te_0\|} \biggr).
\end{eqnarray*}
We deduce that, for any $S$ such that
$ \llVert\pi_{S_0\setminus S}\te_0\rrVert\leq s_n$, denoting by
$c_{\al}$ the normalizing
constant of the density $x\to c_{\al}\exp(-|x|^\alpha)$,
\begin{eqnarray*}
\frac{c_{\al}^{p_n}}{c_{\al}^{|S\cap S_0|}}\frac{\sup_{\q\in
A}g_{S\cap S_0}(\q)} {
\inf_{\q\in B}g_{S_0}(\q)} &\leq&\exp\biggl[\|\te_0
\|_\al^\al\biggl\{ \biggl(1-\frac{3r_n}{4\|\te_0\|}
\biggr)^\al- \biggl(1-\frac{r_n}{4\|\te_0\|} \biggr)^\al
\biggr\} \biggr]
\\
&\leq& \exp\biggl[-2\al(5/8)^{\alpha-1}r_n\frac{\|\te_0\|_\al^\al
}{4\|\te_0\|}
\biggr]\leq\exp\bigl[-4\al(5/8)^{\alpha-1}r_n^2
\bigr],
\end{eqnarray*}
where to obtain the second last inequality we have used that
for any $0\le t\le1/8$ and $\alpha\ge1$
it holds $(1- t)^\alpha-(1-3t)^\al=\int_1^3\alpha t(1-ut)^{\alpha
-1} \,du
\ge2\alpha t (1-3/8)^{\alpha-1}$. Hence the expression in
(\ref{EqQuotientProofLB}) is bounded above by
\begin{eqnarray*}
&&
\sum_S \frac{\Pi_n(S)}{\Pi_n(S_0)}(c_{\al}s_n)^{|S\cap S_0|-p_n}
\frac{v_{|S\cap S_0|}}{v_{p_n}} e^{-4\alpha(5/8)^{\alpha-1} r_n^2}
\\
&&\qquad\leq e^{-4\alpha(5/8)^{\alpha-1} r_n^2}\frac{
e^{Cp_n}}{\Pi_n(S_0)}\sum_{S}
\Pi_n(S)
\\
&&\qquad\le e^{-4\alpha(5/8)^{\alpha-1} r_n^2} e^{Cp_n}e^{cp_n\log(n/p_n)}
\end{eqnarray*}
by Lemma~\ref{lemvol}.
The right-hand side is of smaller order than $e^{-r_n^2}$.
An application of Lemma~\ref{lemlb} concludes the proof for the
case that $1<\alpha<2$.

The proof in the case that $\alpha\ge2$ follows the same lines,
except that we use the inequality $\|\q\|_\alpha\le\|\q\|$, for
every $\q\in\RR^p$, without the factor $p^{1/\alpha-1/2}$ that is necessary
if $\alpha<2$. We define $s_n=(r_n/8)\|\q_0\|_\alpha/\|\q_0\|$.
\end{pf*}

\section*{Acknowledgment}

The authors would like to thank Subhashis Ghosal for suggesting a
simplified argument in the proof of Proposition~\ref{PropositionDimension}.

\begin{supplement}%[id=suppA]
\stitle{Supplement to ``Needles and Straw in a Haystack:
Posterior concentration for possibly sparse sequences''}
\slink[doi]{10.1214/12-AOS1029SUPP} %[doi,text={...}] - jei reikia
%suskaldyti doi
\sdatatype{.pdf}
\sfilename{aos1029\_supp.pdf}
\sdescription{This supplementary file contains the proofs of some
technical results appearing in the paper.}
\end{supplement}

% imsref loaded by lrinkeviciute, 2012-09-19 08:43:05
% imsref loaded by lrinkeviciute, 2012-09-19 08:45:50

\printaddresses

\end{document}